\documentclass[11pt,reqno,a4paper]{amsart}
\usepackage{amssymb}
\usepackage{graphicx}

\def\PSL{\mbox{$\mbox{\rm PSL}_2(\mathbb{Z})$}}
\def\G{\mbox{$\boldsymbol{\mathfrak{G}}$}}
\def\R{\mbox{$\boldsymbol{\mathcal{R}}$}}
\def\Gtilde{\mbox{$\boldsymbol{\widetilde{\mathfrak{G}}}$}}
\def\Ghat{\mbox{$\boldsymbol{\widehat{\mathfrak{G}}}$}}

\def\GAMMA{\mbox{$\boldsymbol{\Gamma}$}}
\def\X{\mbox{$\boldsymbol{X}$}}
\def\Xtilde{\mbox{$\widetilde{\phantom{X}}\!\!\!\!\!\!\boldsymbol{X}$}}

\def\P{\mbox{$\boldsymbol{\mathcal{P}}$}}
\def\ds{\displaystyle}

\newcommand{\mfrac}[2] 
{\raisebox{0.045em}{\mbox{\footnotesize$\displaystyle
\frac{#1}{#2}$}}}

\newcommand{\reff}[1]{\mbox{\rm (\ref{#1})}}

\newcommand{\sss}{\scriptscriptstyle}

\begin{document}
\noindent {\footnotesize {\sc Moscow Mathematical Journal}} \hfill
\texttt{\small http:/$\!$/arXiv.org/math.CA/0205261}\hfill

\noindent {\tiny  Volume 8, Number 2, April--June 2008, Pages 233--271}
\hfill\\\\

\title{On the uniformization of algebraic curves}
\author{Yu.\,V. Brezhnev}
\email{brezhnev@mail.ru}

\hyphenation{uni-for-mi-za-tion}
\subjclass{Primary 30F10; 30F35}


\keywords{Uniformization of algebraic curves, Riemann surfaces,
Fuchsian equations/groups,  monodromy groups,
accessory parameters, modular equations,
conformal maps, curvelinear polygons,
$\theta$-functions, Abelian integrals,
metrics of Poincar\'e, moduli spaces.}

\begin{abstract}
Based on Burnside's parametrization of the algebraic curve
$y^2=x^5-x$ we obtain  remaining attributes of its uniformization:
associated Fuchsian  equations and their solutions, accessory
parameters, monodromies, conformal maps, fundamental polygons, etc.
As a generalization, we propose a way of uniformization of arbitrary
curves by zero genus groups. In the hyperelliptic case all the
objects of the theory are explicitly described. We consider a large
number of examples and, briefly, applications: Abelian integrals,
metrics of Poincar\'e, differential equations of the Jacobi--Chazy
and Picard--Fuchs type, and others.
\end{abstract}

\maketitle
\section{Introduction}
\noindent
As it is understood in contemporary language,
the uniformization of the plane algebraic curves, defined
by  irreducible over $\mathbb C$ polynomial
\begin{equation}\label{1}
\X\!:\;\; F(x,y)=0,
\end{equation}
is, in a broad sense, a 50/50 union of analytic theory and
theory of infinite discontinuous groups of  linear transformations
of the form
\begin{equation}\label{mobius}
\tau\mapsto\frac{a\,\tau+b}{c\,\tau+d}\;.
\end{equation}

({\bf A}) Algebraic and geometric aspects of the problem are rather
well developed. These include theory of Fuchsian and Kleinian groups
\cite{511,lehner,466,355, katok}, their subgroups \cite{466, 468,
483,  P82}, hyperbolic geometry, and fundamental polygons
\cite{lehner, 511, 355, 551}. Descriptions of concrete class of
curves in terms of these objects were considered in the classical
literature and related to the names of Schottky \cite{405}, Weber
\cite{417}, Burnside \cite{303}, and Whittaker \cite{418}. Geometry
of polygons is most studied in the cases of low genera
\cite{367,354,317,buser,aigon} and Fuchsian groups of higher genera
have received  less attention for the reason of computation
complexity. However, pure algebraic and geometric approaches do not
touch the question of explicit parametrization of the implicit
relation (1). This point is the central problem in the
uniformization (see emphasis on pp.\,176--177 in \cite{weyl}) and
has essentially analytical characterization.

({\bf B}) The analytical description implies the presence of two
analytic functions, {\em single-valued\/} in a domain of their existence,
of the {\em global\/} uniformizing parameter $\tau$
\begin{equation}\label{2}
x=\chi(\tau),\qquad  y=\varphi(\tau)
\end{equation}
such that, being substituted into  (1), the functions
\reff{2} turn this equation into identity in the variable $\tau$.
These functions are transcendental and possess the property
of being invariant under  certain  subgroups of the
M\"obius group \reff{mobius} \cite{weyl,lehner}.

First examples of explicit parametrizations of higher genera curves
were obtained by Jacobi in his {\em Fundamenta Nova} \cite{jacobi}.
After Jacobi, his followers (Schl\ae{}fli, Sohnke et all) expanded
this list and, nowadays, it is known as classical modular equations
\cite{483,weber}. Modern achievements in  the explicit construction
of the uniformizing map are also related to particular examples of
modular curves \cite{knopp,mckean,kuyk,atkin}. Besides these names
special mention should be made of non well-known results of the
1900's which deserve to be better known and concern, in the spirit
of Fricke \& Klein \cite[\protect{\bf I}]{483},  to subgroups of the
full modular group \PSL: dissertation by Morris \cite{447}, a series
of papers by Hutchinson \cite{411/2, 411,412,413}, Young
\cite{young1, young2}, and Dalaker \cite{dalaker}. Some of them
already  dealt with modular groups of Jacobi's $\vartheta$-constants
and representations in terms of $\vartheta$- or general
$\boldsymbol{\vartheta}$-constants. The Klein curve
$x^3y+y^3z+z^3x=0$ \cite{klein2} is the famous and most considerable
result  along these lines.

Group of invariancy of the functions \reff{2} is  not unique. It
defines a kind of uniformization. We shall restrict our
consideration only to Fuchsian groups of first kind \cite{511},
because any algebraic curve can be uniformized by such a kind of
groups \cite{449}. We call this group (more precisely its matrix
\mbox{$2\!\times\!2$} representation) as {\em  group $\G$
uniformizing the equation \reff{1}\/}. Without entering upon
details, related to possible punctures and compactification, $\G$ is
one of exact matrix representations of the fundamental group
$\pi_1^{}(\R)$ of a Riemann surface $\R$ of genus $g$ \cite{449}.
The group \G\ is finitely-generated and must have $2\,g$ generators
at the least. Two automorphic functions \reff{2} are generators of
the function field $\boldsymbol{F}=\mathbb{C}(x,y)$
--- field of meromorphic
functions on the curve $\X$ or, which is the same,
on \R\ defined by the equation (1) \cite{lehner}.

\subsection{Outline of the problem}
There are two kinds of the analytical descriptions of uniformizing
functions: 1) those which express the functions \reff{2} in terms of
substitutions of a group and 2) those which express the functions in
terms of coefficients of equation of a curve. $\Theta$-series of
Poincar\'e \cite{P82,303} and the theory of automorphic forms
\cite{551,468,466,467,479,golubev,473,lehner,496,494} solve the
first problem, however the structural form of these series,
presently, does not appear to furnish an effective realization of
the solution. Due to essential non-commutativity of the groups, the
parametric representations for these summations and their
differential properties are unknown. There is no also known way how
to solve these tasks. The second kind is a description in terms of
Fuchsian differential equations  (see last sentence in the book
\cite{551}). The  functions \reff{2} satisfy nonlinear autonomous
{\sc ode}'s of third order
\begin{equation}\label{dex}
\ds[x,\,\tau]=Q_1^{}(x,y), \qquad \ds[y,\,\tau]=Q_2^{}(x,y),
\end{equation}
where $Q_{1,2}^{}$ are certain rational functions of the indicated
variables and $[x,\,\tau]$  is a notation for the meromorphic derivative%
\footnote{\ Right-bracket notation goes back to Dedekind
\cite{dedekind}. See \cite{br} for properties of this object.}
$$
[x,\,\tau] \equiv \frac{1}{x_\tau^2}\,\{x,\,\tau\}\;.
$$
The symbol $\{x,\,\tau\}$, as always, stands for Schwarz's
derivative \cite{golubev,551,P84,395}.
It is well-known  that the equations (\ref{dex}) are equivalent to linear
{\sc ode}'s of Fuchsian type
\begin{equation}\label{Q12}
\Psi_{\mathit{xx}}=\mfrac12 \,Q_1^{}(x,y)\,\Psi, \qquad
\widetilde\Psi_{\mathit{yy}}=\mfrac12
\,Q_2^{}(x,y)\,\widetilde\Psi.
\end{equation}
Explicit formula connection between these equations and (\ref{dex})
is also known. The global uniformizing parameter $\tau$ is a ratio
\begin{equation}\label{ratio}
\tau=\frac{\Psi_2(x)}{\Psi_1(x)}
\end{equation}
of two fundamental solutions $\Psi_{\!1,2}(x)$ of the equation
(\ref{Q12}):
$$
\Psi_1(x)=\sqrt{\mathstrut x_\tau}, \qquad
\Psi_2(x)=\Psi_1(x)\!\int\!\!\mfrac{dx}{\Psi_1^2(x)}=
\tau\sqrt{\mathstrut x_\tau}\;.
$$
{\em The problem of construction of the analytic function
$\chi(\tau)$ is thus a problem of formula representation for
inversion of the ratio \reff{ratio}\/}.

In the context of uniformization, first nontrivial Fuchsian
equations arose in an explicit form in the dissertation by Shottky
\cite[pp.\,57--58]{405} for genus 2 and 3. Since the automorphic
$\Theta$-series (fonctions Fuchsiennes by Poincar\'e) are ones over
the monodromical groups of certain Fuchsian {\sc ode}'s, the central
problem of the theory consists in joining of both the components
({\bf A}) and ({\bf B}). In the framework of differential equations
(\ref{dex}--\ref{Q12}) this means  determination the  rational
function $Q_1^{}(x,y)$ as  function of  equation of a curve. One
part of coefficients is determined algebraically by the equation
\reff{1} itself \cite[\S V]{P84}. To find out the remaining
coefficients (linear combination of holomorphic quadratic
differentials), so that corresponding monodromy group of the
equation shall be Fuchsian, is the celebrated problem of accessory
parameters \cite{P84, smirnov,417, 418, 342, 343}. We call them
$A$-parameters.

There is  a one-to-one correspondence between geometric
characteristics of the fundamental polygon, identification of its
sides, and $3g-3$ $A$-parameters (Klein--Poincar\'e). This number is
equal to dimension of the moduli space $\boldsymbol{\mathcal M}_g$
of algebraic curves \cite{449,bers} and therefore knowledge of the
parameters could give us information on a coordinate representation
of this space. Among other things such representations would
exemplify us nontrivial ($g>1$) explicitly solvable equations of the
Picard--Fuchs type on automorphic forms corresponding to these
groups. Nice examples for genus zero are known  \cite{takh2} but far
from  being exhausted. Very lucid exposition of this subject and
further references can be found in \cite{mckay2,harnad}.

\subsection{State of the art}
Being marked in the literature more than once as ``notorious''%
\footnote{D.\,Hejhal remarks in the survey \cite{343}: {\em
``Although there was extensive work \ldots\ the over-all state of
affairs is still rather primitive $($after 70 some years$)$''\/}.}
\cite{342,hempel,tah5,tah4,315}, the problem of the parameters
($A$-problem), in the analytical sense, is an essentially global one
and transcendental in its general statement
\cite{golubev,P84,467,496,342,343}. There exist however instances of
its exact solution. First  examples were obtained by Sir Edmund
Whittaker \cite{420} and his pupils \cite{421, 318,378}. This list
was  enlarged by R.\,Rankin in the fifties \cite{393, 395}. After
the 1950's no new examples, related to higher genera curves, were
found.

For the elliptic curves, which are  isomorphic to the one canonical
form $y^2=4\,x^3-g_2^{}\,x-g_3^{}$, the solution of the problem is
given by the classical theory of Jacobi--Weierstrass. Namely,
Abelian integrals and meromorphic functions are described by
Weierstrassian functions $\sigma,\zeta,\wp,\wp'$ and the moduli
space --- by Jacobi's $\vartheta$-constants.

Concerning higher genera, despite the fact that there are numerous
examples of modular parametrizations, these are not mentioned in
works on the $A$-problem. On the other hand, these examples would
give a large family of nontrivial solvable Fuchsian equations, but,
as far as we know, none of them has been considered in the spirit of
explicit construction of objects constituting the uniformization in
a straightforward setting: Fuchsian equations, their monodromies,
and Abelian integrals. In this context, the question about
integrals, in fact,  has not been explicitly arisen in the
literature (see however \cite{knopp2}). Meantime this is, in
essence,  the central point in the problem for the simple reason
that the set of meromorphic functions is a particular case of
Abelian integrals and construction of the formers can be replaced by
construction of the meromorphic and logarithmic integrals through
the holomorphic ones \cite{baker}.

In 1893 Burnside found explicit parametrization of a nice algebraic
relation \cite{308}
\begin{equation}\label{bc}
y^2=x^5-x\;.
\end{equation}
In 1958, Rankin, considering Whittaker's approach to
uniformization of hyperelliptic curves \cite{418,420},
found the $A$-parameters for the curve \reff{bc} in the
framework of
this kind of uniformization but uniformizing functions for
Weber--Whittaker- or Rankin's Fuchsian
equations are unknown hitherto.
Burnside's parametrization is not a Whittaker's one.
They have differed Fuchsian equations, but $A$-parameters
have turned out to be coincided and equal to zero.
It was shown in  \cite{br}. In the same place,
we started to develop an analytical description
for the Burnside parametrization
and discussed a relationship to a Whittaker's conjecture
\cite{420,421}. Rankin, in one of his last papers \cite{397},
considered group
properties of Burnside's parametrization and
presently,  the example \reff{bc} is the only one
having all the attributes of the theory  in an explicit form. This is
one of the subjects of the present work.

The following table exhibits  the state of the art for
the cases $g=0,1$ and $g>1$.

\medskip

\begin{center}
\begin{tabular}{l|c|c} \hline \hline
\multicolumn{1}{c|}{${}_{\ds\mathstrut}^{\ds\mathstrut}$\small\sc
problem}&
\small$g=0,\;1$ & \small $g>1$\\
\hline
\small 1. \it Explicit parametrization${}_{}^{\ds\mathstrut}$ & $+$ & $-$ \\
\small 2. \it Determining Fuchsian {\sc ode}'s
and their solutions& $+$ & $-$ \\
\small 3. \it Poincare's and $q$-series& $+$ & $\pm$\\
\small 4. \it Group-geometric properties & $+$ & $\pm$ \\
\small 5. \it Conformal maps and
polygons&$\pm$&$-$ \\
\small 6. \it Uniformization  in
the language of $\Theta$-functions & $+$ & $-$\\
\small 7. \it  Abelian integrals as functions of the parameter $\tau$ &+&$-$\\
\small 8. \it Covers of higher genera curves\,$:$
{\sc ode}'s, $\Theta$'s, etc & & $-$\\
\small 9. \it Moduli spaces and applications${}_{\mathstrut}^{}$ &
\;\;\,\,\,\,$+$\,\small$(!)$ &
\;\;\,\,\,\,$-$\,\small$(!)$\\
\hline\hline
\end{tabular}
\end{center}

\medskip

We displayed not all existing problems because the theory has
numerous branches and direct links to diverse areas. Pertinent
discussions and additional bibliography can be found in the surveys
\cite{468,343,bers}, monographs \cite{449,539,473} or classical
books \cite{551,klein3,466,467,479,golubev,483,451,lehner}. It
should be emphasized here that all the above-mentioned points, being
developed for the elliptic curves, provide enormous number of their
applications. For the case of $g>1$, general proofs of
existence-theorems, based on the concept of a covering surface
(simply connected version of \X), were completed in the 1900's by
Koebe and Poincar\'e.

\subsection{Sketch of the work}
Content  of the work is new and known results are accompanied by
corresponding references.

Sect.\,2 contains important results concerning transformations
between Fuchsian equations,  definition of automorphisms of elements
of the field  $\boldsymbol{F}$ as matrix monodromies and other
terminology. Sects.\,3--5 deal with the Burnside parametrization of
the curve \reff{bc} for purposes of subsequent generalizations. One
constructs explicitly bases of all monodromies, Fuchsian equations,
their solutions, group properties, conformal mappings, fundamental
polygons,   and distributions of exact values of uniformizing
functions therein (inversion problem). New kind of uniformization of
the hyperelliptic curves by free groups is described in sect.\,6. We
construct corresponding Fuchsian equations and bases of their
monodromies. The next sect.\,7 is an illustration of close relation
between hyper- and non-hyperelliptic curves in examples by
$\vartheta$-constant uniformizations and modular equations. Sect.\,8
concludes all the preceding stuff with the main theorem about
hyperelliptic and zero genus characterization of proposed
uniformization. Consequences of the main theorem and extra examples
are considered in sect.\,9. The final sect.\,10 contains comments to
applications with explicit formulas following from the
considerations above and brief historic-bibliographical remarks.
{\sc Appendix} contains some known and new identities between
Jacobi's $\vartheta$-constants. To avoid lengthening the reference
list we refrain from wider discussions of classical and modern
literature in this field. This will be a subject matter of a
detailed survey.

\section{Transformations between Fuchsian equations}

\noindent Since $\chi(\tau)$ and $\varphi(\tau)$ are connected by
the algebraic equation \reff{1}, the {\sc ode}'s
(\ref{dex}--\ref{Q12}) are not independent of each other. The
present section clarifies such  connections.

\subsection{Lemma} {\em The functions $Q_{1,2}^{}(x,y)$  are
connected by the linear relation}
\begin{equation}\label{q12}
\big(Q_1^{}+\{F,\,x\}\big)\mfrac{F_y}{F_x}-
\big(Q_2^{}+\{F,\,y\}\big)\mfrac{F_x}{F_y}=
3\Big(\!\ln\!\mfrac{F_x}{F_y}\Big)_{\!\mathit{xy}}\quad
\mbox{mod}\;\;\big\langle F(x,y)\big\rangle.
\end{equation}

\noindent {\it Proof.} Differentiating the identity
$F\big(\chi(\tau),\varphi(\tau)\big)=0 $ three times with respect
to $\tau$ and making use of  equations (\ref{dex}) we get three algebraic
relations on $\tau$-derivatives.
Elimination the  derivatives of the functions $x,y$ leads to
the formula \reff{q12}. \hfill \rule{0.5em}{0.5em}

In terms of  Fuchsian equations  the {\bf Lemma} is a particular
case $z=y$ of invariant  form of the general transformation law of
the meromorphic derivative under the change $x \mapsto z$ \cite{br}:
\begin{equation}\label{law}
\left\{
\begin{array}{l}
\ds\widetilde\Psi_{\mathit{zz}}=\mfrac12\Big\{[z,\,x]
+\mfrac{1}{z_x^2}\,Q_1^{}(x,y)\Big\}\widetilde\Psi,\qquad
\widetilde\Psi(z)=\sqrt{\mathstrut z_x}\,\,\Psi(x)\\\\
\ds[z,\,\tau] = [z,\,x]+\frac{1}{z_x^2}\,Q_1^{}(x,y)
\end{array}\right.,
\end{equation}
therefore solvability of one of the equations \reff{Q12} involves
solvability of the other.

Thus the $A$-problem is a problem of {\em one\/} Fuchsian equation:
{\em no matter to which function the Fuchsian equation
corresponds\/}. It may be one of the generators or  arbitrary
element $z=R(x,y)$ of the field $\boldsymbol{F}$. Clearly the
element $z$ should be chosen to be ``good'' for subsequent
integration of the equations.

Formula \reff{q12} was obtained for an algebraic dependence between
generators $x$ and $y$, but, in form \reff{law}, {\bf Lemma} may be
applied to the arbitrary change of variables. Indeed, the Riemann
surface \R\  is always a ramified cover and there are infinitely
many representations for such covers. These may be algebraic, as
\reff{1}, or transcendental. The last case occurs if the curve (1)
is realized as a cover of a Riemann surface of genus higher than
zero. A large number of examples of covers of elliptic tori is
provided by the theory of elliptic solitons \cite{acta}. In the
language of Fuchsian equations, transitions between these
$(1\!\mapsto\! N)$-covers correspond to  the choice of a generator
of the field $\boldsymbol{F}$ and are merely the changes \reff{law}
of independent variable $x\mapsto z$ and linear transformation
$\Psi(x)\mapsto\widetilde\Psi(z)$. The transcendental
representations correspond to non-algebraic changes and, thereby,
generate Fuchsian equations  on surfaces of higher genera. First
example of solvable Fuchsian equation on a torus was exhibited in
\cite[\S\,4]{br}. We shall refer to the {\bf Lemma} in the form
\reff{law}.

\subsection{Automorphisms of meromorphic functions on \R}
Let $\G_x$ and $\G_y$ be
matrix monodromy groups (in short, the monodromies)
of the equations \reff{Q12} respectively.
Assuming the $A$-parameters to be known, $\G_x$ and $\G_y$
become automorphism groups  of the functions $\chi(\tau)$ and
$\varphi(\tau)$
$$
\G_x\equiv\mbox{{\bf Aut}}\big(\chi(\tau)\big), \qquad
\G_y\equiv\mbox{{\bf Aut}}\big(\varphi(\tau)\big)
$$
by definition. Since domains of existence of analytic functions
$\chi(\tau)$ and $\varphi(\tau)$ coincide (upper $(\tau)$-half-plane
$\mathbb{H}$ under certain normalization), the matrix group $\G$ is
the intersection $\G_x \bigcap \G_y$ being non-elementary Fuchsian
group of genus $g$. The groups $\G_x$ and $\G_y$ do not necessarily
coincide but may be conjugated by a non-degenerate matrix
$\G_y=U\G_xU^{-1}$, be subgroups of one another, and have differed
genera including genus zero. In other words, the uniformizing group
\G, monodromies of generators $\G_{x,y}$ or some element of the
field $\G_z$ are not the same. Conversely, every Fuchsian group
$\G_x$, represented geometrically by some polygon, may be considered
as a monodromical group of certain Fuchsian equation in $x$, i.\,e.
as an automorphism group of  the analytic function $x=\chi(\tau)$.
This is direct consequence of a positive solvability of the
Riemann--Hilbert problem for {\sc ode}'s of 2nd order: given a
monodromy there exists the corresponding Fuchsian equation
\cite[\S\,IV]{P84}.

\noindent{\bf Corollary 1.}
{\em All representations for the
$A$-parameters are computed using  \reff{law}\/}.

It is not at first sight apparent but later, in sect.\,8, we shall
see that for arbitrary algebraic curve, under suitable
uniformization, there exists a generator  with a zero genus
automorphism $\G_x$ and rational function $Q(x)$. In other words,
the usual  approach to the problem of uniformization through the
analysis of a Fuchsian equation with {\em canonical\/} system of
$2g$ generators corresponding to a monodromy $\G_x=\G$ of genus $g$
is noneffective and rather inconvenient way. There always exists
uniformization by a group of genus zero. In any case, $\R$, $\G$,
and a fundamental polygon $\P$ are constructed by minimal number of
copies of $\P_{\!x}$ or $\P_{\!y}$ in their own planes. This number
is the number of sheets of \R.

The transformations described above do not solve the problem of
searching for the global monodromies  or even monodromy of one
equation if  the monodromy of the other is known\footnote{Fuchsian
{\sc ode} itself has no  problem of parameters. These are defined by
an external condition: what properties of  $\Psi$-function would we
like to have? For example: $(\bullet)$ $\Psi$ is an algebraic
function (Riemann (1850's), Schwarz, Klein, Brioschi, Jordan
(1870's)) or, more generally; $(\bullet)$ $\Psi$ is representable by
indefinite quadratures (Liouville--Picard--Vessiot theory);
$(\bullet)$ the {\sc ode} \reff{Q12} brings about a conformal
mapping of the upper half-plane onto  given a shape circular polygon
or, which is our case; $(\bullet)$ when the monodromy is infinite
and Fuchsian?}. These transformations are important for reducibility
of the problem to simpler one. For example for obtaining larger
groups. A larger group is easier described since corresponds to a
simpler Fuchsian equation or even, if possible, to equation for a
triangle group. The first case, where it was done, is an example of
Whittaker \cite{420} and its subsequent generalizations \cite{421,
318, 378, 379}. It was found, that both groups $\G_{x,y}$ have genus
zero and the functions $Q_{1}^{},\,Q_{2}^{}$ are rational functions
of $x,\,y$ respectively. Moreover, $\G_y$ is a 2-generated triangle
group and $\G$ is of index $2g+1$ subgroup of $\G_y$. This point was
a crucial step which allowed to compute global monodromies of  the
equations \reff{Q12} for the curves $y^2=x^{2g+1}+1$. Accessory
polynomials have turned out to be equal to zero. See \cite{br} for
additional comments about these reductions. Below we shall consider
the equation (\ref{bc}) in detail, generalize it, and establish its
analogy with Whittaker's approach, including the hypergeometric
reducibility.

\section{Monodromy $\G_x$  for Burnside's parametrization}

\noindent Structure of  the group  $\G_x$ of the function
$\chi(\tau)$ was recognized and used by Burnside in work \cite{308}.
It was found \cite{308, 397} that $\G_x$  is
the principal congruence subgroup $\GAMMA(4)$ of level 4 in the
full modular group $\GAMMA(1)=\PSL$ \cite{494}. It is known
\cite[{\bf I}: p.\,355]{483} that $\GAMMA(4)$ is generated by the
five matrices
\begin{equation}\label{g4}
\mbox{$\ds V_{0\ldots 4}^{}=\left\{ \left(\!\!
\begin{array}{rl}
1 & 4\\
0 & 1
\end{array}\!\!
\right), \; \left(\!\!
\begin{array}{rl}
\phantom{-}1 & 0\\
-4 & 1
\end{array}\!\!
\right), \;
 \left(\!\!
\begin{array}{rl}
3 & -4\\
4 & -5
\end{array}\!\!
\right), \;
\left(\!\!
\begin{array}{rl}
7 & -16\\
4 & -9
\end{array}\!\!
\right), \;
 \left(\!\!
\begin{array}{rl}
11 & -36\\
4 & -13
\end{array}\!\!
\right)\right\}$}.
\end{equation}
Fuchsian equation \reff{Q12} for the function $x=\chi(\tau)$ has the
following form \cite{br}
\begin{equation}\label{fex}
\Psi_{\mathit{xx}}=-\frac
14\,\frac{x^8+14\,x^4+1}{(x^5-x)^2}\,\Psi
\end{equation}
and its monodromy $\G_x$ is $\GAMMA(4)$. Index
$\big|\GAMMA(1)\!:\!\GAMMA(4)\big|=24$ follows from  degree of the
algebraic relation between $\chi(\tau)$ and the absolute modular
invariant $J(\tau)$ \cite{br}
\begin{equation}\label{j}
J(\tau)=\frac 1{108} \frac{(x^8+14\,x^4+1)^3} {(x^5-x)^4}.
\end{equation}
Substitution (\ref{j}) into the equation (\ref{fex}), being
considered as algebraic one $x\mapsto J$, and simultaneous changing
$\Psi(x) \mapsto \widetilde \Psi(J)=\sqrt{J_x\mathstrut}\,\Psi(x)$
({\bf Lemma}) reduce the equation (\ref{fex}) to the well-known
equation (Bruns 1875)
$$
\widetilde{\Psi}_{\!\mathit{JJ}}^{}=
-\frac{36\,J^2-41\,J+32}{144\,J^2\,(J-1)^2}\,\widetilde{\Psi}
$$
and vice versa. Another change of variables $\big(x,\Psi(x)\big)
\mapsto \big(z, \Phi(z)\big)$:
\begin{equation}\label{eq8}
z=x^4, \qquad \Psi(x)=y\,\Phi(z)
\end{equation}
reduces the equation (\ref{fex}) to the hypergeometric equation of
Legendre \cite{PL,golubev,watson}
\begin{equation}\label{eq9}
z(z-1)\,\Phi_{\mathit{zz}}+(2z-1)\,\Phi_z+\mfrac 14\,\Phi=0\,.
\end{equation}
Clearly, the function $z(\tau)$ is any
of Legendre's moduli:
\begin{equation}\label{k2}
[z,\,\tau ]=-\frac12\,\frac{z^2-z+1}{z^2 (z-1)^2} \quad
\Rightarrow\quad z=k'{}^2\;\mbox{or}\;
\;k^2\Big(\mfrac{a\tau+b}{c\tau+d}\Big).
\end{equation}
Solutions in terms of the hypergeometric function ${}_2F_1$
have the following form
$$
\Psi_1(x)=\sqrt{x^5-x\,} \;{}_2F_1\!\!\left(\mfrac12, \mfrac
12;1\Big| x^4\right),\qquad \Psi_2(x)=\sqrt{x^5-x\,}
\;{}_2F_1\!\!\left(\mfrac 12, \mfrac 12;1\Big| 1-x^4\right)
$$
and can be written down in terms of complete elliptic integrals $K$
and $K'$:
\begin{equation}\label{psi}
\begin{array}{l}
\ds\Psi_1(x)=\frac{2}{\pi}\,\sqrt{x^5-x\,}\,\int\limits_{\!0}^{\,1}
\!\!\mfrac{d\lambda}{\sqrt{(1-\lambda^2) (1-x^4\lambda^2)}}=\qquad\,
\frac{2}{\pi}\,\sqrt{x^5-x\,}\;K(x^2),\\
\ds\Psi_2(x)=\frac{2}{\pi}\,
\sqrt{x^5-x\,}\,\int\limits_{\!0}^{\,1^{\ds\mathstrut}}
\!\!\mfrac{d\lambda}{\sqrt{(1-\lambda^2) (1-(1-x^4)\lambda^2)}}=
\frac{2}{\pi}\,\sqrt{x^5-x\,}\;K'(x^2).
\end{array}
\end{equation}
The integral form of solutions to Fuchsian equations is
convenient for obtaining  corresponding
monodromy groups \cite{golubev,PL}. On the other hand,
$\tau$-representations
for the $\Psi$-functions are important as explicit examples of
single-valued formulas on \R.
The integrals $\Phi$, as functions
of the ratio $\tau$,  read as \cite{jacobi,watson}
\begin{equation}\label{Ktheta}
\Phi_1^{}(\tau)=\frac{\pi}{2}\,\vartheta_3^2(\tau), \qquad
\Phi_2^{}(\tau)=-i\,\tau\,\frac{\pi}{2}\,\vartheta_3^2(\tau)\,.
\end{equation}
See sect.\,10.4 for the $\tau$-representation of the functions
$\Psi_{1,2}$. It is well-known (the result by L.\,Fuchs; see also
\cite[pp.\,130--132]{496}, \cite{PL,golubev}) that the monodromy
group of equation (\ref{eq9}) is $\GAMMA(2)$ and is generated by the
two matrices $\ds\big(\begin{smallmatrix}1 & 0_{}\\2 &
1\end{smallmatrix}\big)$, $\ds\big(\begin{smallmatrix}1 & 2_{}\\0 &
1\end{smallmatrix}\big)$. The group $\G_x$ is of  index 4 subgroup
of $\GAMMA(2)$ which follows from degree of the substitution
(\ref{eq8}) or a  property of $\GAMMA(4)$ of being a subgroup of
free group $\GAMMA(2)$.

The equation (\ref{fex}) has five finite regular singularities $
e_{1\ldots5}^{}=\big\{0,\,\pm 1,\,\pm i\big\} $ and a singularity
at infinity $e_6^{}=\infty$. Locally, each singularity corresponds
to one of the generators $S_k$ of $\G_x$. Obviously $S_1\cdots
S_5=S_\infty^{-1}$ and there are only five independent
generators. As a set of group
generators, $\big\{S_k\big\}$ is equivalent to the matrices
(\ref{g4}).

Fundamental polygon for $\GAMMA(4)$ is a 10-gon $\P_{\!x}$ with
all parabolic vertices ({\sc Fig.\,1}). It has
genus zero \cite{483,lehner} and the function $\chi(\tau)$ has
one simple pole therein. A natural question arises:
what function does realize
a conformal mapping of  $\P_{\!x}$ with a proper correspondence
of these borders? Section 5 answers this question.

The next problem is to determine the group $\G$ uniformizing the
equation (\ref{bc}).

\begin{figure}[htbp]
\centering
\includegraphics[width=8cm]{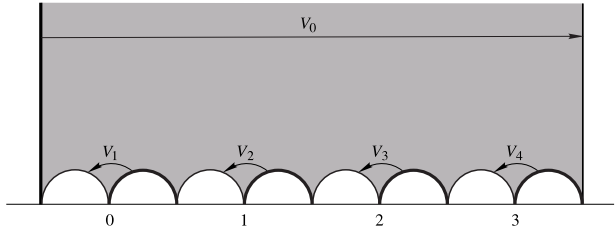}
\caption{Fundamental 10-gon $\P_{\!x}$ for the group
$\G_x=\GAMMA(4)$.}
\end{figure}

\section{Group $\G$ for Burnside's parametrization}
\noindent The Riemann surface $\R$ of the equation (\ref{bc})
consists of two $x$-leaves and, the  fundamental polygon $\P$
for $\G$ must contain two arbitrary neighbouring copies of
$\P_{\!x}$.

\noindent
{\bf Proposition.} {\em The group \G\ of the equation
$(\ref{bc})$ is a free subgroup of rank\/ $9$ and
index $2$ in $\GAMMA(4)$. It is generated by the
following elements\,$:$
{\small $$
\begin{array}{c}
T_0^{}=V_0\,V_0^{}=\left(\!\!
\begin{array}{ll}
1 & 8\\
0 & 1
\end{array}\!\!
\right)
\\
\begin{array}{ll}
T_1^{}=V_0^{}\,V_1^{}=\left(\!\!
\begin{array}{rl}
15 & -4\\
4 & -1
\end{array}\!\!
\right)^{\ds\mathstrut}_{\ds\mathstrut} \qquad &
T_5^{}=V_0^{}\,V_1^{-1}=\left(\!\!
\begin{array}{rl}
17 & 4\\
4 & 1
\end{array}\!\!
\right)
\ \\
T_2^{}=V_0^{}\,V_2^{}=\left(\!\!
\begin{array}{rl}
19 & -24\\
4 & -5
\end{array}\!\!
\right)^{\ds\mathstrut}_{\ds\mathstrut} &
T_6^{}=V_0^{}\,V_2^{-1}=\left(\!\!
\begin{array}{rl}
21 & -16\\
4 & -3
\end{array}\!\!
\right)
\
\\
T_3^{}=V_0^{}\,V_3^{}=\left(\!\!
\begin{array}{rl}
23 & -52\\
4 & -9
\end{array}\!\!
\right)
&
T_7^{}=V_0^{}\,V_3^{-1}=\left(\!\!
\begin{array}{rl}
25 & -44\\
4 & -7
\end{array}\!\!
\right)^{\ds\mathstrut}_{\ds\mathstrut} \
\\
T_4^{}=V_0^{}\,V_4^{}=\left(\!\!
\begin{array}{rl}
27 & -88\\
4 & -13
\end{array}\!\!
\right)^{\ds\mathstrut} & T_8^{}=V_0^{}\,V_4^{-1}=\left(\!\!
\begin{array}{rl}
29 & -80\\
4 & -11
\end{array}\!\!
\right).
\end{array}
\end{array}
$$}}%

Since the function $y=\varphi(\tau)$ is a
modular function belonging to $\G_x$ with a multiplier system
$\varkappa=-1$ \cite{494,397}, the automorphism group $\G_y$
coincides with the group $\G$\footnote{Generally, hyperelliptic
$\G_y$ is not equal to  $\G$ as Whittaker's examples show
\cite{420,421,378,318}.}. We expect to see this in the language of
differential equations.

\noindent {\em Proof\/}.
Following the {\bf Lemma}, one gets an equation with
algebraic coefficients
\begin{equation}\label{eq13}
\widetilde\Psi_{\!\mathit{yy}}(y)= -\mfrac14\,
\frac{y^{16}-192\!\cdot\! 5^{-3}\,x\,y^{14}+576\! \cdot\!
5^{-4}\,x^2\,y^{12}+ \cdots +
2^{16}5^{-10}}{y^2(y^8-2^85^{-5})^2}\,\widetilde\Psi(y).
\end{equation}
The equation (\ref{eq13}) has 9 finite singular points and a
singularity at infinity
$$
y=\big\{0, \;2\!\sqrt[-8]{5^5},\;\infty \big\}
$$
which are, except for $y=0$,  branch points of the cover $y \mapsto
x$. Since the point $y=0$ is a usual point (no permutation of sheets
at this place) with  indicial equation $n(n-1)+\frac14=0$, it
produces a  parabolic generator $T_0$ of the {\em global\/}
monodromy $\G_y$. Each of remaining points corresponds to
permutations of sheets $(x_j^{},y) \mapsto (x_k^{},y)$: two of five
$y$-sheets for $E_k=2\!\sqrt[-8]{5^5}$ and all five ones for
$y=\infty$. Indeed, the equation \reff{eq13} has the form (we reduce
an accessory part that does not affect the local analysis)
$$
\widetilde\Psi_{\mathit{yy}}=\left\{\!- \frac14\frac{1}{y^2}
+\frac{3}{2^{10}}\,{\sum}_E\! \mfrac{2^5 x^4-50E^4x^2+ 5^3 E^6
x-2^85^{-1}}{(y-E)^2}-\frac{3\!\cdot\!5^5}{2^{13}}\, {\sum}_E
\!\mfrac{E^7 x^4-\cdots}{y-E} \right\}\widetilde\Psi
$$
and indicial equation for $y=\infty$ takes the form $n^2+n=-\frac14$
with one multiply root $\left\{-\frac12\right\}$. In a
neighbourhood of the points $y=E_k$ the equation takes the form
$$
\widetilde\Psi_{\mathit{yy}}=
\frac{\alpha+\cdots}{(y-E_k)^2}\,\widetilde\Psi,
$$
where  parameter $\alpha$ is equal to, depending on a sheet, one
of the values $\left\{0,0,0,-\frac{3}{16} \right\}$. This yields
an identical transformation of three sheets and
permutation of two other ones mentioned above.
One-fold closed paths around $E_k$ do not generate  elements of a
local monodromy. Nevertheless the latter  is constructed from
these permutations. Accordingly, hyperbolic generators
of the global $\G_y$
can not be recognized by a local analysis but must exist for a
group of genus $g>1$. See \cite{thome} and later works by Thom\'e
for systematic discussion of equations having algebraic
coefficients.

The parabolic generator $T_0$ entails  punctures on
$\R=\mathbb{H}/\G_y$ at five points $y=0$ and makes the group
to be free. Hence  one of ten possible generators should
be defined by other ones: $\G_y$ is a free group of rank 9.

Straightforward calculation of the global monodromy $\G_y$
is difficult enough. At present, in fact, there is no
known method how  to do this for equations with
algebraic (and even rational) coefficients by  means of a computational procedure.
However we may use  two-sheetness of the cover (\ref{bc}). A simplest
way to get  index 2, if group is free,  is to divide the group
into elements of odd and even length with respect to any base. We
take a subgroup with elements of even length over the base
\reff{g4}:
\begin{equation}\label{base}
\G:\qquad V_j\,V_k, \quad V_j\,V_k{}^{\!\!-1}, \quad
V_j{}^{\!\!-1}\,V_k, \quad V_j{}^{\!\!-1}\,V_k{}^{\!\!-1}, \qquad
(j,\,k=0 \ldots 4).
\end{equation}
Perhaps, Rankin tacitly used this fact in \cite{397} to describe
some group properties of the Burnside parametrization. Each matrix
of the group \G\ is equal  to one of the four matrices
$\ds\big(\begin{smallmatrix}1&0_{}\\0&1\end{smallmatrix}\big)$,
$\ds\big(\begin{smallmatrix}1&4_{}\\4&1\end{smallmatrix}\big)$,
$\ds\big(\begin{smallmatrix}5&4_{}\\0&5\end{smallmatrix}\big)$,
$\ds\big(\begin{smallmatrix}5&0_{}\\4&5\end{smallmatrix}\big)$
modulo 8 \cite[{\bf I}: p.\,652]{483}. Clearly, one hundred matrices
\reff{base} is not a minimal base. Since $\GAMMA(4)$ is free, we
obtain,  by the Nielsen--Schreier formula \cite{LS,494}, that
subgroup of index $m$ in a free group of rank $k$ is generated by
$N=(k-1)\,m+1$ elements. Rank of \G\  is equal to $(5-1)\,2+1=9$ as
stated above.

The parabolic generator $T_0$ must double the {\sc Fig.\,1}. This
also suggests how  to construct  unknown generators. Let $\P$ be a
polygon that consists of two copies of $\P_{\!x}$:
$\P=\P_{\!x}\bigcup V_0^{}(\P_{\!x})$. We build the set of 9
elements
\begin{equation}\label{T}
\langle\boldsymbol T\rangle= \big\langle V_0{}^{\!2},\,\,V_0\,V_k,
\,\,V_0\,V_k{}^{\!\!-1}\big\rangle, \qquad k=1 \ldots 4.
\end{equation}
These identify the sides of $\P$ as it is shown in {\sc Fig.\,2}.
\begin{figure}[htbp]
\centering
\includegraphics[width=11cm]{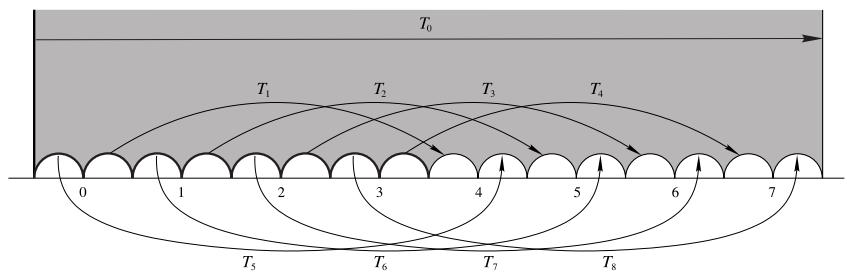}
\caption{$\G$ and $\R$: fundamental 18-gon $\P$ for Burnside's
parametrization.}
\end{figure}
The hyperelliptic involution
$\varphi(V_k\tau)=-\varphi(\tau)$ is obvious and we finish the proof:
$\G=\G_y=\langle \mbox{\textbf{\textit{T}}}\rangle$.
\hfill\rule{0.5em}{0.5em}

We note that the Reidemeister--Schreier rewriting algorithm
\cite{LS} provides another general recipe for obtaining a minimal
set of generators but does not provide a direct geometric
information similarly to {\sc Fig.\,2}. Topologically, all the nine
generators $\{ T_k\}$ (one parabolic and eight hyperbolic)
correspond to five punctures on $\R$ plus canonical base of four
cycles $(\boldsymbol{\mathfrak
{a}}_k^{},\boldsymbol{\mathfrak{b}}_k^{})$ existing on any
compactified Riemann surface of genus 2. Sixth puncture kills the
standard relation $\prod [\boldsymbol{\mathfrak
{a}}_k^{},\boldsymbol{\mathfrak{b}}_k^{}]=1$ between
$\boldsymbol{\mathfrak {a}}$ and $\boldsymbol{\mathfrak{b}}$'s.

Notice that the local analysis described above is partially
simplified if   to use solution
$\widetilde\Psi=\sqrt{\varphi_\tau^{}}$ itself as a function on
$\R$ (which is, however, being constructed):
$$
\widetilde\Psi(\tau)^2=
\frac{5\,\chi(\tau)^4-1}{2\,\varphi(\tau)}\,\chi_\tau^{}(\tau).
$$
Its singularities correspond to $\varphi_\tau^{}=\{0,\infty\}$ and,
as earlier, lead to  the points $y=0$, infinity, and eight points
$5\,x^4-1=0$. Globally, on the corresponding \R's through the
variable $\tau$, all the Fuchsian equations are indistinguishable
and become ({\bf Lemma}) the equation $\Psi_{\tau\tau}=0\cdot\Psi$
with the solution $\Psi=A\tau+B$. The monodromies, in this view, are
also indistinguishable since we fall into the universal cover.

\section{Conformal maps}
\subsection{Burnside's function as a conformation}
The functions $\chi(\tau)$ and $\varphi(\tau)$, as
meromorphic functions on $\R$, assume every complex value the same
number of times: $\chi(\tau)$
--- twice and $\varphi(\tau)$ --- five times.
Owing to the analytic properties, the function $\chi(\tau)$ takes  all
possible values from $\overline{\mathbb C}$ one time for all $\tau
\in \P_{\!x}$. This corresponds to analytical points $(x,y)$. The
second copy of $\P_{\!x}$ corresponds to the analytical points
$(x,-y)$. Since the polygon $\P_{\!x}$ is a simply connected
region, its image $\chi(\P_{\!x})$ must coincide with  the extended
plane $(x)$. The conformality is lost at
$e$-points\footnote{$A$-point of analytic function is a point
where the function takes  the value $A$.} of the function
$\chi$ where $\chi'(\tau)=0$. On the other hand, all sides of
$\P_{\!x}$ are separated to equivalent pairs. Hence  images of
two equivalent sides coincide in the plane $(x)$ and form
a cut between $e$'s. This gives  walking there and back direction
along the cut between $e_k^{}$ which are critical points of the
conformal mapping.  Thus, starting from some point $x$ on a
cut we shall come back to $x$ passing subsequently branch points
$e_{k+1}^{},\,e_{k+2}^{},\,\ldots$. The possible kinds of the cuts
in the plane $(x)$  satisfying the demands above are the star- or
tree-like ones with branches, their degenerations, etc.
One $N$-folded straight line corresponds to
the well-known half-plane constructions of Klein--Schwarz.
From the point of view
of Fuchsian {\sc ode}'s, one star-like pictures are simplest because
generators of the local  monodromy, being only elliptic or parabolic,
become the
same for a global monodromy. It is readily recognized from the
local one. No  hyperbolic operators appear in a set generating $\G_x$
since only neighbouring
sides in the $\tau$-plane are coupled. Burnside's
example  happens to be the very case:

{\em The Burnside's function $\chi(\tau)$ realizes  one-to-one
conformal mapping of the interior $\P_{\!x}$ in $\tau$-plane onto
the exterior of a star-like domain in the plane $(x)$. The center
of the star is the point $e_1^{}=1$  corresponding to the equivalent
$\tau$-points $
\left\{-\frac{1}{2},\,\frac12,\,\frac32,\,\frac52,\,\frac72\right\}$
$($see {\sc Fig.\,1}$)$  forming one vertex cycle of
$\,\P_{\!x}$}.

A shape of the cuts can not be found using  group properties or
defining {\sc ode}, but it is  completely determined by  analytic
function itself. In our case we have \cite{308,br}
\begin{equation}\label{burnside}
\chi(\tau)=\frac{\wp(1|2,2\tau)-\wp(2|2,2\tau)}
{\wp(\tau|2,2\tau)-\wp(2|2,2\tau)}.
\end{equation}
The shape depends on a boundary of fundamental polygon $\P_{\!x}$
which is not unique. In {\sc Fig.\,3}, with a correct scale of arcs,
we use more illustrative disc model obtained from the upper
half-plane of the variable $\tau$ ({\sc Fig.\,2}) by the
transformation $\tau \mapsto \tilde\tau$:
$$
\tilde\tau=i\,\mfrac{\tau-1-i}{\tau-1+i}.
$$

\begin{figure}[htbp]
\centering
\includegraphics[width=11cm]{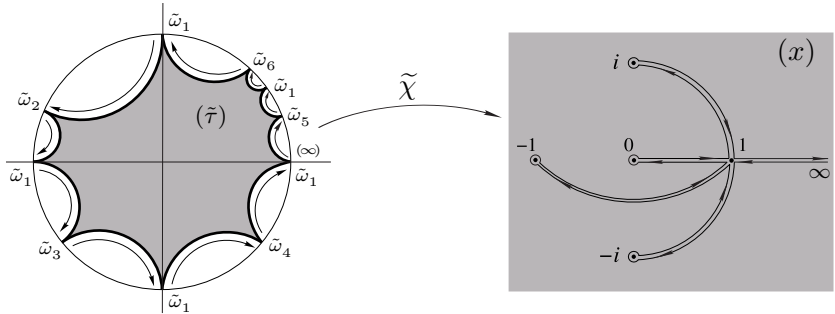}
\caption{The function $\widetilde\chi(\tilde\tau)$ as a
conformation.}
\end{figure}

\noindent The function
$$
\chi\!\left(\!\mbox{\small$(1-i)$}\mfrac{\tilde\tau+1}{\tilde\tau-i}\right)=
\widetilde\chi(\tilde\tau)
$$
brings about a conformal mapping shown in {\sc Fig.\,3}. A big arc
$(-1\ldots 1)$ is of radius $\sqrt{2}$ arc. As one would expect, all
the cuts  have a circular shape and illustrate geodesics between
Weierstrass points in $(\tau)$- and $(x)$-planes (see also
sect.\,10.6).

\subsection{The problem of inversion} Problem of determination of
points on the curve \reff{bc} is equivalent to solution of the
transcendental equation $\chi(\tau_{\sss 0})=A$. Except few points,
the solution $\tau_{\sss 0}$ is unique in $\P_x$. Determination of
$y$ is trivial due to {\bf Proposition}. Making use of Burnside's
function \reff{phi} we have that $\pm y=\{\varphi(\tau_{\sss 0}),
\varphi(\tau_{\sss 0}+4)\}$. The remaining calculations are similar
to those in the elliptic modular inversion problem \cite{watson}.
Taking into account a multi-valuedness of the integrals \reff{psi},
we should choose  the functions $\Psi_{1,2}(x)$  so that their ratio
gets into desired domain independently of the value $x$. We take
$\mathbb{H}$ and compute the value $\tau'$ modulo $\GAMMA(1)$ by the
formula
$$
\tau'= i\,\frac{K(A^2)}{K'(A^2)}
\qquad
\left(=i\,\mfrac{{}_2F_1\!\big(\frac12,\frac12;1;A^4
\big)}{{}_2F_1\!\big(\frac12,\frac12;1;1-A^4 \big)}\right),
\qquad \tau' \in \mathbb{H}
$$
whereupon the sought for point $\tau_{\sss 0}$ is chosen from 24
values $\tau_{1\ldots24}^{}$ being images of $\tau'$ in
$\boldsymbol{\mathcal{P}}_x$ under $\GAMMA(1)$. This generates
the octahedral irrationality $\GAMMA(1)/\GAMMA(4)$ \cite{308,483,klein3}:
\begin{equation}\label{J}
J\big(\tau_{1\ldots24}^{}\big)=\frac{1}{108}\,
\frac{(A^8+14A^4+1)^3}{(A^4-1)^4A^4}\;.
\end{equation}
Exact values of the function $\chi(\tau)$ in
some  points $\tau\in \P_x$ are obtained with the
help of algebraic relation \reff{J} and known  values of
Klein's invariant. Some of these points suggest themselves. These are
vertices of the  $\GAMMA(1)$-triangle
$\tau=\big\{\infty,i, \frac{-1+i\sqrt{3}}{2}\big\}$
corresponding  to the branch points $x=\{0,\,\pm1,\,\pm i,\,\infty\}$
and the values
$$
x_{1\ldots12}^{}=\big\{\!\pm1\pm\sqrt{2\,},\;\pm\sqrt{\pm i\,},\;
\pm i\pm i\sqrt{2\,}\big\}\;,\qquad
x_{1\ldots8}^{}= \sqrt[\leftroot{-1}4]{\!-7\pm
4\,\sqrt{3\,}}
$$
respectively.
Other exact values of $\chi(\tau)$ can be found by involving modular relations
between $J(\tau)$ and $J(n\tau)$, i.\,e. cases where $J$ is an
algebraic number. For example, twenty four $\GAMMA(1)$-images of the point
$\tau=i\sqrt{2}$ generate roots of the equation
$$
(A^8+14A^4+1)^3=500\,(A^4-1)^4A^4.
$$
In particular, $\chi\big(i\sqrt{2}\big)=\sqrt{\!\sqrt{2}-1}$. More
examples are obtained by making use of tables in books on imaginary
quadratic number fields, e.\,g. \cite[{\bf III}]{weber},
\cite{mckean}.

In order to find an inverse image of the real axis of the plane
$(x)$ we use the fact that $J(\tau)$ assumes all possible real
values on a boundary of the well-known fundamental quadrangle for
the group $\GAMMA(1)$. Hence, as it follows from \reff{j}, the
function $\chi(\tau)$ takes all possible real values on arcs forming
a decomposition of $\P_{\!x}$ into  triangles of the group
$\GAMMA(1)$. Analyzing such a decomposition, one can find a curve
$\gamma$ in the plane $(\tilde\tau)$ such that $\chi(\gamma)$ is a
real quantity. It is shown in {\sc Fig.\,4}.

In addition, in {\sc Fig.\,4}, we pictured incomplete triangulation
of one sheet $\boldsymbol{\mathcal{P}}_x$ of the Riemann surface
$\R$ and  images of some its $n$-gons on the plane $(x)$: 3-gon $A$,
5-gon $B$, and the shaded 6-gon. As a corollary we get an inverse
image of the upper half-plane of the variable $x$.
\begin{figure}[htbp]
\centering
\includegraphics[width=11cm]{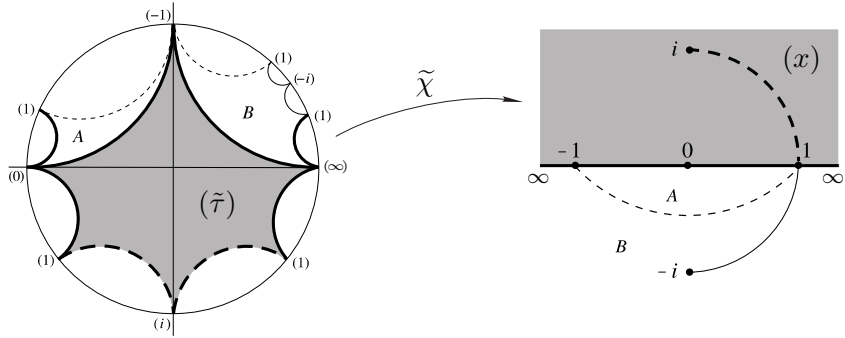}
\caption{Conformal mapping onto upper half-plane $(x)$.}
\end{figure}

\section{Uniformization of hyperelliptic curves}

\subsection{Brief comments to known uniformizations}
Uniformizations of Schottky \cite{405} and Weber \cite{417}
deal with
multi-connected domains and this point
complicates  analysis  of
corresponding differential equations. It also
entails that holomorphic and meromorphic Abelian integrals are
multi-valued functions of the parameter $\tau$.

Uniformization of Whittaker \cite{418} is free from such  features
and is the most natural generalization of $\wp$-parametrization of
elliptic curves. Simply connected  polygon for the group
$\G_x$ has a minimal number of sides
$4g$, the group contains no elliptic or parabolic elements, but
is not free and has the only defining relation
\begin{equation}\label{whitrel}
T_1\,T_2\,\cdots\,T_{2g}= T_{2g}\,\cdots\,T_2\,T_1,
\end{equation}
where $T_k$ are properly sorted generators identifying opposite
sides.

An essential advantage of Weber's and Whittaker's uniformizations
is that the Fuchsian linear {\sc ode} has
rational coefficients in  the variable $x$. See also
\cite[p.\,57]{405} for genus 2.  It should be noted here that all
known  examples of solved Fuchsian equations \cite{420,318,393}
correspond to the
rational functions $Q_1^{}(x)$.

In the next section
we shall show that these facts are not accidental
and there always exists a zero genus uniformization with a
simply connected  polygon.
Presence of punctures makes a group to be free
in a complete analogy with Burnside's example.

We emphasize  that not all curves admit  uniformization by a
conformation of upper half-plane and  subsequent application of
Schwarz's reflection principle likewise the celebrated 14-gonal
picture of Klein \cite{klein3, haskell} or Weber's case \cite{417}.
Meanwhile every uniformizing function or any algebraic function
$R(x,y)$, being  meromorphic  on $\R$, implements a conformal map
onto a whole plane $\overline{\mathbb{C}}$ with corresponding cuts
between critical points of its own mapping. Nontrivial
counterexamples, when the reflection principle may not be applied,
arise even in the elliptic case \cite[{\bf IV}:
pp.\,143--166]{tannery}. The {\sc Fig.\,4} shows that Burnside's
example is another one along these lines.

\subsection{Parabolic uniformization}
Let us
take $2\,g+2$ arbitrary real points $\omega_k^{}$ such that
$\omega_1^{}<\omega_2^{}<\cdots<\omega_{2g+2}^{}$. Let these
points, except for $\omega_1^{}$, be fixed points of $2\,g+1$
parabolic transformations
$\big(V_1,\,V_2,\,\ldots,\,V_{2g},\,V_0\big)$  and
$\varepsilon_k^{}$ be $2g$ arbitrary real points located between
the points $\omega_k^{}$ as follows

\bigskip
\centerline{
\unitlength=1mm
\begin{picture}(100,17)
\put(-2,5){\line(1,0){100}}
\put(2,5){\circle*{1.}}\put(1,1){\footnotesize$\omega_1^{}$}
\put(13,5){\circle*{1.}}
\put(12,1){\footnotesize$\omega_2^{}$}
\put(11.8,9){\footnotesize$V_1^{}$}
\put(22,5){\circle*{0.6}}\put(21,1){\footnotesize$\varepsilon_1^{}$}
\put(32,5){\circle*{1.}}
\put(31,1){\footnotesize$\omega_3^{}$}
\put(30.8,9){\footnotesize$V_2^{}$}
\put(41,5){\circle*{0.6}}\put(40,1){\footnotesize$\varepsilon_2^{}$}
\put(48,5){\circle*{1.}}
\put(47,1){\footnotesize$\omega_4^{}$}
\put(46.8,9){\footnotesize$V_3^{}$}
\put(54,0.5){$\cdots\cdots$}
\put(54,9){$\cdots\cdots$}
\put(70,5){\circle*{1.}}
\put(68,1){\footnotesize$\omega_{2g+1}^{}$}
\put(68.3,9){\footnotesize$V_{2g}^{}$}
\put(82,5){\circle*{0.6}}
\put(80,1){\footnotesize$\varepsilon_{2g}^{}$}
\put(91,5){\circle*{1.}}
\put(89,1){\footnotesize$\omega_{2g+2}^{}$}
\put(89,9){\footnotesize$V_0^{}$}
\put(82,16){\footnotesize$(\tau)$}
\end{picture}}

\noindent We construct the parabolic transformations $V_k$ identifying
arcs of circles being orthogonal to the real axis by the
rule ($\omega_{2g+2}^{}$ may be thought of as $\infty$  in {\sc Fig.\,1}):
\begin{equation}\label{ident}
\begin{array}{l}
\stackrel{V_1}
{\stackrel{\!\!\!\!\!\scalebox{4.5}[1.5]{\mbox{{\large
$\ds\curvearrowleft$}}}}{\big[\omega_1^{},\,\omega_2^{}\big]\;
\big[\omega_2^{},\,\varepsilon_1^{}\big]}_{{}_{\ds\mathstrut}}},

\quad \stackrel{V_2}
{\stackrel{\!\!\!\!\!\scalebox{4.5}[1.5]{\mbox{{\large
$\ds\curvearrowleft$}}}}{\big[\varepsilon_1^{},\,\omega_3^{}\big]\;
\big[\omega_3^{},\,\varepsilon_2^{}\big]}_{{}_{\ds\mathstrut}}},
\quad\ldots\ldots\\
\ldots\ldots,\quad

\stackrel{\;\;\;\;V_{2g}}
{\stackrel{\;\scalebox{4.5}[1.5]{\mbox{{\large
$\ds\curvearrowleft$}}}}{\big[\varepsilon_{2g-1}^{},\,\omega_{2g+1}^{}\big]\;
\big[\omega_{2g+1}^{},\,\varepsilon_{2g}^{}\big]}_{{}_{\ds\mathstrut}}},
\quad

\stackrel{\;V_{0}}
{\stackrel{\!\!\!\scalebox{4.5}[1.5]{\mbox{{\large
$\ds\curvearrowleft$}}}}{\big[\varepsilon_{2g}^{},\,\omega_{2g+2}^{}\big]\;
\big[\omega_{2g+2}^{},\,\omega_1^{}\big]}_{{}_{\ds\mathstrut}}}.
\end{array}
\end{equation}

\noindent Since we deal only with Fuchsian groups,  parabolic
transformation contains two real undetermined constants. What
number of the parameters do we have at our disposal? It is equal
to $2\,g+1$ real fixed points of generators $V_k$ plus $2\,g$
arbitrary real points $\varepsilon_k^{}$ which define second
parameters of generators $V_k$. One remaining point $\omega_1^{}$
is equivalent to the point $\varepsilon_{2g}^{}$ by the
construction:
$$
V_0\,\big(V_1\,V_2\cdots V_{2g}\big)\,\varepsilon_{2g}^{}=
V_0\,\omega_1^{}=\varepsilon_{2g}^{}.
$$
Hence this scalar identity kills one of $2\,(2\,g+1)$ parameters in
the set $\{V_k\}$. All total is $4\,g+1$. The group in question is
not essentially different from one which is obtained by
transforming it with arbitrary real transformation. Subtracting 3
real parameters, one gets $4\,g-2$  real parameters, that is $2\,g-1$
complex parameters, and we relate them to branch points of a
hyperelliptic curve
\begin{equation}\label{hyper}
y^2=(x-e_1^{})\cdots(x-e_{2g+1}^{})
\end{equation}
after the normalization $e_1^{}=0,\;e_2^{}=1$ and
$e_{2g+2}^{}=\infty$. We set $\G_x=\big\langle V_k \big\rangle$.

Furthermore, we construct  $x=\chi(\tau)$ as a conformal mapping of
interior of the $(4\,g+2)$-side polygon in the plane $(\tau)$ onto
the plane $(x)$ with star-like cuts (see {\sc Fig.\,3}). It is clear
that points $\omega$ and $\varepsilon$ are critical points of the
mapping. On the other hand, the meromorphic derivative $[\chi,\tau]$
approaches  infinity only at critical points of  analytic functions:
$\chi'(\tau)=0$ \cite{br}. Therefore, as $(\omega,\varepsilon)$'s
lie on the real axis, we make the well-known change of a local
parameter \cite{551}
$$
\tau\mapsto q:\qquad \ln q=\frac{2\,\pi\, i}{\tau-\omega}
$$
and, repeating arguments of the works \cite{418,br}, one gets
\big($x(\tau)=X(q)$\big)
$$
[X,\,q]-\frac{1}{2\,q^2}\frac{1}{X_q^2}
=-\frac12\,\frac{1}{(X-e_k^{})^2}+\cdots.
$$
The equation, satisfied by the  function
$x=\chi(\tau)$, thus acquires the form
\begin{equation}\label{general}
[x,\,\tau]= -\frac12 \left\{\sum_{k=1}^{2g+1}\frac{1}{(x-e_k^{})^2}-
\frac{2\,g\,x^{2g-1}+A(x)}{(x-e_1^{})\cdots(x-e_{2g+1}^{})}\right\},
\end{equation}
where $A(x)=a_1^{}x^{2g-2}+\cdots+a_{2g-1}^{}$ is the accessory
$A$-polynomial. The local uniformizing parameter $\tau$
(quasi-uniformising variable in Whittaker's terminology
\cite{418})  becomes the global one under an appropriate
choice of the polynomial $A(x)$.

Now we construct  group $\G$ of the equation \reff{hyper}. We
generate a second copy of $(4\,g+2)$-gon $\P_{\!x}$ by the
transformation $T_0^{}=V_0^2$ and erase the boundary which separates
$\P_{\!x}$ and its copy. The identification of sides of built
$(8\,g+2)$-gon $\P=\P_{\!x}\bigcup V_0(\P_{\!x})$ is done as in
\reff{T}:
$$
\G:\qquad\big\langle\boldsymbol T\big\rangle= \big\langle
V_0{}^{\!2},\,\,V_0\,V_k, \,\,V_0\,V_k{}^{\!\!-1}\big\rangle,
\qquad k=1 \ldots 2g.
$$
The last step is to count the Eulerian characteristics $ E=c-s+1$
wherein $c$ is a number of cycles, $2\,s$ is number of sides of
$\P$, and $1$ is a number of polygons. We have $2\,s=(4\,g+2)\,2-2$.
All the points $\omega_{\sss 2\ldots 2g+2}$ form $2\,g+1$  parabolic
cycles and $(\omega_1^{},\,
\varepsilon_1^{},\,\ldots,\,\varepsilon_{2g}^{})$ form one long
cycle: $c=(2\,g+1)+1$. Checking
$$
E=\big((2\,g+1)+1\big)- \mfrac{(4g+2)\,2-2}{2}+1=2-2\,\mathit{genus}
$$
completes the construction: {\em genus\/} is equal to  $g$.

Note that the case  $g=2$ is the only one that corresponds to maximal
possible number of sides of a normal polygon $2\,s=12\, g+4\!\cdot\!
0-6=18$ \cite[p.\,174]{weyl}, \cite{511}.
Thus, proposed uniformization can be thought of as a parabolic
version of Whittaker's one
\begin{equation}\label{whittaker}
[x,\,\tau]= -\frac38 \left\{\sum_{k=1}^{2g+1}\frac{1}{(x-e_k^{})^2}-
\frac{2\,g\,x^{2g-1}+A(x)}
{(x-e_1^{})\cdots(x-e_{2g+1}^{})}\right\},
\end{equation}
but without defining relation \reff{whitrel}. Indeed, the polygons
$\P_{\!x}$ have the same number of sides $(4\,g+2)$ for these
uniformizations and, as some examples show, the conjecture of
Whittaker (when it does hold) is applicable to this version also.
This observation arises a question: {\em when the $A$-polynomials
coincide for the parabolic equation \reff{general} and Whittaker's
one \reff{whittaker}?} See \cite{br} and recent paper \cite{girondo}
on Whittaker's conjecture.

\section{Relation between hyper- and non-hyperelliptic curves. Examples.}

\noindent In this section we shall exhibit that one function may
uniformize many curves by means of subgroups of its own automorphism
group. Let us consider Jacobi's curve
\begin{equation}\label{jacobi}
z^4+w^4=1.
\end{equation}
Based on different representations for the
function \reff{burnside} \cite{397,br}, results of sect.\,3, and
$\vartheta$-constant identities \cite{tannery,weber,watson} we conclude
that
\begin{equation}\label{1997}
\chi(\tau)= -\frac{\vartheta_4(\tau/2)}{\vartheta_3(\tau/2)}
=-\sqrt{k'}\Big(\mfrac{\tau}{2}\Big).
\end{equation}
On the other hand, Legendre's moduli  $k(\tau)$ and $k'(\tau)$
$$
k(\tau)=\frac{\vartheta_2^2(\tau)}{\vartheta_3^2(\tau)},\qquad
k'(\tau)=\frac{\vartheta_4^2(\tau)}{\vartheta_3^2(\tau)}
$$
are single valued functions and arbitrary algebraic/logarithmic
functions of them with arbitrary ramifications only at the points
$\{k,\,k'\}=\{0,\,\pm 1,\,\infty\}$, say
\begin{equation}\label{kk}
\sqrt[\leftroot{0}\uproot{1}4]{k'(\tau)}=
\frac{\vartheta_4^{}(2\tau)}{\vartheta_3^{}(\tau)},\quad
\sqrt[\leftroot{0}\uproot{1}6]{kk'(\tau)}=
\sqrt[\leftroot{0}\uproot{1}6]{4}\,
\frac{\eta(\tau)}{\vartheta_3^{}(\tau)},\quad
\sqrt[\leftroot{-2}\uproot{2}n]{k^p(k-1)^m(k'+1)^s(\tau)},\quad\ldots,
\end{equation}
will be single-valued also. The only question is their formula
representation.  In particular, for Burnside's function $y=\varphi(\tau)$
\cite{308} (see \cite{br} for notations and details)
{\small
\begin{equation}\label{phi}
\varphi(\tau)=4i\,\mfrac {\big[\wp(\tau)-\wp(2\tau)\big]\!
\big[\wp\big(\frac{\tau}{2}\big)-\wp(\tau)\big]\!
\big[\wp\big(\frac{\tau}{2}\big) -\wp(\tau+2)\big]\!
\big[\wp\big(\frac12\big)-\wp(2\tau+1)\big]\!
\big[\wp\big(\frac12\big)-\wp(1)\big]_{{\mathstrut}}}
{\big[\wp\big(\frac \tau 2\big)-\wp(1)\big]^{{\mathstrut}}\!
\big[\wp\big(\frac \tau 2\big)-\wp(2\tau+1)\big]\,\wp'\big(\frac
12\big)\, \wp'(\tau)}
\end{equation}}
we have
\begin{equation}\label{y}
y=\sqrt{x(x^4-1)}=
\sqrt[\leftroot{1}\uproot{1}4]{k'}\!\left(\mfrac{\tau}{2}\right)
k\!\left(\mfrac{\tau}{2}\right)^{\ds\mathstrut}.
\end{equation}
The formula \reff{1997} was obtained in \cite{397}. We make the
change $(x,\,y,\,\tau) \mapsto (-x, \,-i\,y,\,2\tau)$, as always in
the sequel, whereupon we arrive at  pending ``root-free" formulas.

{\em The function of Burnside
\reff{burnside} realizes uniformization of both curves \reff{bc},
\reff{jacobi} and explicit parametrizations
have the following form
\begin{equation}\label{prop2}
\left\{ x=\frac{\vartheta_4(\tau)}{\vartheta_3(\tau)},\quad
y=4\,i\,\frac{\eta^3(2\tau)}{\vartheta_3^3(\tau)}
\right\}, \qquad \left\{
z=\frac{\vartheta_4(\tau)}{\vartheta_3(\tau)},\quad
w=\frac{\vartheta_2(\tau)}{\vartheta_3(\tau)}\right\}.
\end{equation}
The $\eta,\vartheta$-quotient for the function $y$ in \reff{prop2}
satisfies the second equation in \reff{dex} with $Q_2^{}$ defined
by \reff{eq13} and the quotients $(x,\,z,\,w)$ satisfy one
equation
\begin{equation}\label{z4}
[x,\,\tau]= -\frac12\frac{x^8+14x^4+1}{x^2(x^4-1)^2}
\end{equation}
with automorphism group $\G_x = \G_z = \G_w=\GAMMA(4)$.
Corresponding solutions $\Psi_{1,2}$ of Fuchsian equation \reff{Q12}
under normalization
$\Psi_1^2(x)=\frac{d}{d\tau}\frac{\vartheta_4(\tau)}{\vartheta_3(\tau)}$
have the form
\begin{equation}\label{Psi}
\Psi_1(x)={\textstyle\sqrt{\mfrac{i}{\pi}}}\,\sqrt{x^5-x}\,K'(x^2),
\qquad
\Psi_2(x)=i\,{\textstyle\sqrt{\mfrac{i}{\pi}}}\,\sqrt{x^5-x}\,K(x^2).
\end{equation}}

\noindent {\em Proof\/}. Formula for
the variable $y$ is obtained  by combining  \reff{y},
first equality in \reff{kk},  doubling formulas \reff{vartheta},
and  a relation between Dedekind's function and $\vartheta$-constants:
$2\,\eta^3(\tau)=\vartheta_2(\tau)\vartheta_3(\tau)\vartheta_4(\tau)$.
Isomorphism of monodromies  follows from
the well-known modular transformations for the $\vartheta$-constants
\cite{tannery,watson}
$$
w(\tau)=z\Big(\mfrac{-1}{\tau}\Big)=x\Big(\mfrac{-1}{\tau}\Big)\;.
$$
The equations (\ref{z4}--\ref{Psi}) are straightforward consequences of
\reff{fex} and rules of differentiation of the $\vartheta$-constants
\reff{fundamental}. \hfill\rule{0.5em}{0.5em}

As a direct generalization, we uniformize  Fermat's curves
\begin{equation}\label{fermat}
z^n+w^n=1
\end{equation}
in the sense that the generators $(z,w)$ of this field satisfy one Fuchsian equation
\begin{equation}\label{n}
[z,\,\tau]= -\frac12\frac{z^{2n}+(n^2-2)\,z^n+1}{z^2(z^n-1)^2}\;.
\end{equation}
Another generalization is the set of higher curves of Whittaker
$k^2=1-z^n$ (compare with Whittaker's conjecture \cite{420,421,br}).
The function $k(\tau)$ satisfies the Heun's equation
\begin{equation}\label{heun}
[k,\,\tau]= -\frac12\,\frac{(k^2+1)^2}{k^2(k^2-1)^2}\;.
\end{equation}
Explicit solutions to equations \reff{fermat} are provided by
suitable $\vartheta$-quotients for the functions
\mbox{$z=\!\sqrt[\leftroot{-2}\uproot{2}n]{k^2(\tau)},\;
w=\!\sqrt[\leftroot{-2}\uproot{2}n]{{k'}{}^2(\tau)}$} though this is
rather nontrivial problem\footnote{We suggest that the general
$\vartheta$-constants of the form
$\theta\big(f(\tau)\big|g(\tau)\big)$ solve the problem.}.
Associated Fuchsian {\sc ode} \reff{Q12} is solved in terms of
hypergeometric function
${}_2F_1\!\big(\frac12,\frac12;1\big|\sqrt[\uproot{2}n]{z\,}\big)$.
Since all the singularities in \reff{n} are parabolic ones, the
relevant $q$-series, having numerous number-theoretic applications,
can be obtained directly from this equation even if we do not have
$\vartheta$-formulas for them. See \cite{br} for example $n=4$.

Internal points $\Im(\tau)>0$ also admit exact series expansions.
The simplest way to do this is to make use of the algebraic relation
\reff{j} between $x$ and $J$, because the series for Klein's
invariant $J(\tau)$ and Hurwitz's functional $q$-series for the
forms $g_{2,3}^{}(\tau)$ are well studied. For example for the point
$\tau=i$ we have
$$
J(\tau)=1-12\,\mfrac{g_2^{}}{\pi^2}\,(\tau-i)^2-
12\,i\,\mfrac{g_2^{}}{\pi^2}\,(\tau-i)^3+
\mfrac{g_2^{}}{\pi^2}\!\left(\mfrac{184}{3}\mfrac{g_2^{}}{\pi^2}+9
\right)\! (\tau-i)^4+\cdots,
$$
$$
\chi(\tau)=1-\sqrt{2}+a\,(\tau-i)+
\Big(\frac{\mathstrut}{\mathstrut}\!\!
i\,a+\mfrac{3\,\sqrt{2}-4}{2\,\pi^2}\,g_2^{}
\!\!\frac{\mathstrut}{\mathstrut}\Big)(\tau-i)^2+ \cdots\;,
$$
where
$$
a^2=\mfrac{4\,\sqrt{2}-6}{\pi^2}\,g_2^{},\qquad g_2^{} \equiv
g_2^{}(i)=\pi^4\Big\{\mfrac{1}{12}+20\,
\sum\limits_{k=1}^{\infty}\mfrac{k^3}{e^{2\pi k}-1}
 \Big\}=11.817045\ldots\;.
$$
As a byproduct, these series generate different kind  exact
identities of the form
$$
\frac{\vartheta_4^{}}{\vartheta_3^{}}\!\Big(\mfrac i2\Big)=
\sqrt{2}-1,\qquad
\vartheta_4^8(2i)=\mfrac{8}{\pi^4}\,g_2^{}(i),\quad\ldots\;
\mbox{etc}.
$$

The example \reff{jacobi} seems quite obvious and
symmetrical due to Jacobi's quartic identity
$\vartheta_3^4(\tau)=\vartheta_2^4(\tau)+\vartheta_4^4(\tau)$.
The next one may be more appropriate \cite{jacobi,
dedekind,tannery}. It gives solution of the equations
\reff{fermat} and \reff{n} under $n=8$:
\begin{equation}\label{8}
\left\{ \frac{\vartheta_4^{}(4\tau)}{\vartheta_3^{}(2\tau)}
\right\}^8+ \left\{\! \frac{1}{\sqrt{2}\,}\,
\frac{\vartheta_2^{}(\tau)}{\vartheta_3^{}(2\tau)} \right\}^8=1,
\end{equation}
where we have used Landen's transformation in $\vartheta$-form:
$2\,\vartheta_2^{}(2\tau)\,\vartheta_3^{}(2\tau)=\vartheta_2^2(\tau)$
(Cayley calls \reff{8}, in a form of product, the remarkable
identity).

The equations \reff{n} look as ``hyperelliptic'' equation
\reff{general}. Indeed, let $\alpha$'s be branch points in the
equations (\ref{hyper}--\ref{general}). Under $n=4$ the equations
\reff{z4} and \reff{n} become
$$
[z,\,\tau]=
-\frac12\left\{\sum\limits_{j=1}^5\frac{1}{(z-\alpha_j^{})^2}-
\frac{4\,z^3+0}{z^5-z}\right\},\qquad
\alpha_j^{}=\big\{0,\,\sqrt[\uproot{1}4]{1}\big\}
$$
with zero accessory polynomial $A(z)$ for the Burnside curve
$y^2=z^5-z$ \cite{br}.  The same situation  takes place for
higher $n$. Under $n=8$ we have
\begin{equation}\label{4}
[z,\,\tau]=
-\frac12\left\{\sum\limits_{j=1}^9\frac{1}{(z-\alpha_j^{})^2}-
\frac{8\,z^7+0}{z^9-z}\right\}, \qquad
\alpha_j^{}=\big\{0,\,\sqrt[\uproot{1}8]{1}\big\}.
\end{equation}
Monodromy $\Gtilde_z$ of this equation has genus zero.
It uniformizes the hyperelliptic curve
$\Xtilde:\; y^2=z^9-z$. The function $z(\tau)$ is
given by any of terms in \reff{8}, say
\begin{equation}\label{kkk}
z(\tau)=\pm\frac{\vartheta_4^{}(2\tau)}{\vartheta_3^{}(\tau)}.
\end{equation}
(The same function uniformizes the curve $k^2=1-z^8$). The second function
$$
y(\tau)=\sqrt{z(\tau)^9-z(\tau)}= i\,\sqrt{z(\tau)}\,k(\tau)=
i\,\sqrt[\leftroot{1}\uproot{1}8]{k'}(\tau)\,k(\tau)=\cdots
$$
is a single-valued one by construction above. Indeed
$$
\begin{array}{l}
\ds\cdots= i\,\frac{\vartheta_2^2(\tau)}{\vartheta_3^3(\tau)}\,
{\sqrt{\vartheta_4^{}(2\tau)\,\vartheta_3^{}(\tau)}}^{\ds\mathstrut}=i\,
e^{\frac{\pi i}{24}(\tau+1)}_{\mathstrut}\,
\frac{\vartheta_2^2(\tau)}{\vartheta_3^3(\tau)}\, \eta_3^{}(\tau)
\prod\limits_{n=1}^{\infty}\!\! \big(1-q^{4n-2}\big)^{\!\frac12}
\end{array}
$$
and the product $\boldsymbol{\prod}\big(1-q^{4n-2}\big)$ never
vanishes at $\mathbb{H}$. The monodromy $\Gtilde_y$ has genus 4,
index 2 in 9-generated group $\Gtilde_z$ of genus zero, and rank
$N=(9-1)2+1=17$. An explicit representation for this function is a
separate question. The series expansions for $y(\tau)$ are obtained
with use of Fuchsian {\sc ode} \reff{4} or from the last
$\vartheta,\eta$-formula. The group $\Gtilde_y$ itself is
constructed as shown in the previous section. All the groups are
free. Moreover well-known classical Jacobi--Sohnke's modular
equations of arbitrary level $m$  relate moduli
$\sqrt[\leftroot{1}\uproot{1}4]{k}(\tau)$ and
$\sqrt[\leftroot{1}\uproot{1}4]{k}(m\tau)$. They are also
uniformized by the ``hyperelliptic'' function (\ref{4}--\ref{kkk})
of genus zero with just described hyperelliptic curve $\Xtilde$  of
genus $g=4$. For example Jacobi's equations of levels 3 and 5  read
as
\begin{equation}\label{modular4}
z^4-w^4-2zw(1-z^2w^2)=0, \qquad
\left\{z=\frac{\vartheta_4^{}(6\tau)}{\vartheta_3^{}(3\tau)},\qquad
w=\frac{\vartheta_4^{}(2\tau)}{\vartheta_3^{}(\tau)}\right\},
\end{equation}
\
\begin{equation}\label{modular8}
z^6-w^6+5z^2w^2 (z^2-w^2)+4zw(1-z^4w^4)=0,
\end{equation}

$$
\left\{z=\frac{1}{\sqrt{2}}
\frac{\vartheta_2^{}(5\tau)}{\vartheta_3^{}(10\tau)},\qquad
w=\frac{1}{\sqrt{2}}\frac{\vartheta_2^{}(\tau)}{\vartheta_3^{}(2\tau)}
\right\}.
$$
See \cite[{\bf III}]{weber}, \cite{klein3}, and
\cite{mckean,farkas2000,atkin} for additional information about
modular equations. Described non-hyperelliptic examples belong to
Jacobi \cite{jacobi} and, as we have seen now, these all are
directly connected to suitable hyperelliptic curves. All  said above
holds for modular curves of fractional levels
$F\big(\sqrt[\leftroot{1}\uproot{1}4]{k}(n\tau),
\sqrt[\leftroot{1}\uproot{1}4]{k}(m\tau)\big)=0$ or relations
between quantities $\sqrt[\leftroot{1}\uproot{1}4]{k}(m\tau)$,
$\sqrt[\leftroot{1}\uproot{1}2]{k}(m\tau)$, \ldots
$\sqrt[\leftroot{-1}\uproot{2}p]{k}(m\tau)$ which are obtained by
eliminations  common functions
$\sqrt[\leftroot{1}\uproot{1}4]{k}(\tau),\;\ldots$ etc in various
combinations. Corresponding monodromies  will have again zero genus
but they are not  necessarily   conjugated.

Additional identities among $\vartheta$- and $\eta$-constants
generate extra examples of parameterizations. For example two
identities of Dedekind \cite[p.\,191--192]{dedekind}
$$
2^4\,\eta_1^8(\tau)+\eta_2^8(\tau)+ e_{\mathstrut}^{\frac23\pi
i}\,\eta_3^8(\tau)=0\,, \qquad
\eta_1^{}(\tau)\,\eta_2^{}(\tau)\,\eta_3^{}(\tau)= e^{\frac{\pi
i}{24}}_{\mathstrut}\,\eta^3(\tau)
$$
yield, after some simplifications, many non-obvious  hyperelliptic
and other curves ($\eta$-identities on p.\,407 in \cite{farkas2000}
are not true). One of such examples is the curve
\begin{equation}\label{new}
y^2=3\,x^5+10\,x^3+3\,x\,,\qquad \left\{
x=\frac{\vartheta_3^2+2\,\eta^2} {\vartheta_3^2-2\,\eta^2}\,,\quad
y=4\,\frac{(\vartheta_4^4-\vartheta_2^4)\,\vartheta_3^2}
{\big(\vartheta_3^2-2\,\eta^2\big){}^3} \right\}\,.
\end{equation}
We also note that differentiation the known
$\vartheta,\eta$-identities (with the help of rules
\reff{fundamental}) generate new ones.  The literature on modular
equations is enormous
\cite{kuyk,weber,knopp,klein3,lehner,483,far,farkas2000,468,atkin}
and we stop exemplifying. In the next section we shall show that
every algebraic curve can be uniformized in a zero genus manner.

\section{Uniformization of arbitrary curves}
\noindent The possibility   of uniformization of arbitrary curves in
a unified manner was pointed out already by Whittaker \cite{420} in
the framework of groups  generated by self-inverse substitutions
$S_k^2=1$. It is not difficult to see that his approach is
transformable onto groups described in  sect.\,6 with advantages of
relation-free groups. We use a property that every plane algebraic
curve can be birationally transformed into a curve with only double
points. This is implemented with the help of some generic linear
transformation $u'_j=A_{jk}^{}\,u_k^{}$ in homogeneous coordinates
$(u_1^{}:u_2^{}:u_3^{})$.

\subsection{Main theorem}
In order to explain further ideology, we consider the Klein curve
$Y^2(Y-1)=X^7$ as an example. Originally, in ``native'' variables
of $\GAMMA(1)/\GAMMA(7)$, this curve is rather complicated
168-sheet cover. Klein found its compact  form
$u_1^3\,u_2^{}+u_2^3\,u_3^{}+u_3^3\,u_1^{}=0$ \cite{klein1}, which
is  isomorphic to the previous one.

Performing, as an instance, the following birational transformation
$$
\mbox{\footnotesize
$\ds
\left\{x=-\frac{XY+Y}{X^3+XY},\;\;
y=\frac{5X^3+4XY+Y}{X^3+XY}\right\},
\quad
\left\{ X=\frac{x+y-5}{x-y+5}, \;\;
Y=\frac{(x+y-5)^3}{(x+y-3)(x-y+5)^2} \right\}
$}\;,
$$
we get four-leaved cover of $(x)$-plane
\begin{equation}\label{*}
\X:\;\;y^4-4\,(x+6)\,y^3-6\,(x^2-6\,x-33)\, y^2-\cdots{}
{}+12\,x+865=0
\end{equation}
with 12 simple branch-points $\alpha_k^{}$. These points are
roots of the discriminant equation
$$
D(x):\quad 31\,x^{12}+324\,x^{11}+1746\,x^{10}+\cdots -24\,x+23=0.
$$
Let us think of   $\alpha_k^{}$ as branch points of the
hyperelliptic curve
\begin{equation}\label{hyperKlein}
\begin{array}{l}
\ds
\Xtilde:\;\;z^2=31\,x^{12}+324\,x^{11}+1746\,x^{10}+
5916\,x^9+\cdots
-24\,x+23={}_{\ds\mathstrut}^{}\\
\ds
\phantom{\Xtilde:\;\;z^2}=31\,\big(x-\alpha_1^{}\big)\cdots \big(x-\alpha_{12}^{}\big).
{}_{\ds\mathstrut}^{}
\end{array}
\end{equation}
We call the curve $\Xtilde$ the hyperelliptic curve {\em
associated\/} to  given a curve $\X$ (in our case, Klein's curve). The
multi-valued function
$$
\sqrt{31\,x^{12}+324\,x^{11}+1746\,x^{10}+
5916\,x^9+\cdots -24\,x+23\,}
$$
can be made a single-valued one if we shall construct the function
$x=\chi(\tau)$ as a conformation with, at least, 2-folded
$\alpha$-points $x=\alpha_k^{}+c_k^{}(\tau-\tau_k^{})^{2m}+\cdots$
or with a parabolic behaviour of the type
$x=\alpha_k^{}+q+c_k^{}q^2+\cdots$ \cite{br}. Based on the section
6, such a function is readily constructed as satisfying the Fuchsian
{\sc ode}'s \reff{general} or \reff{whittaker}. Modifying the
equation (\ref{general}) for the case of even number of finite
points
\begin{equation}\label{even}
[x,\,\tau]= -\frac12
\left\{\sum_{k=1}^{2g+2}\frac{1}{(x-\alpha_k^{})^2}-
\frac{2\,(g+1)\,x^{2g}+2\,g\,{\sum}\,\alpha_k^{}\!\cdot
x^{2g-1}+A(x)}
{(x-\alpha_1^{})\cdots(x-\alpha_{2g+2}^{})}\right\},
\end{equation}
we get an equation for the function $x$ uniformizing both Klein's curve
\reff{*} and the associated curve  \reff{hyperKlein}
$$
[x,\,\tau]= -\frac12
\left\{\sum_{k=1}^{12}\frac{1}{(x-\alpha_k^{})^2}-
\frac{12\,x^{10}-10\,\frac{324}{31}\, x^9+A_8(x)}
{(x-\alpha_1^{})\cdots(x-\alpha_{12}^{})}\right\}
$$
with some accessory polynomial $A_8(x)$.  The monodromy
$\Gtilde_x$ of this equation (chosen  to be Fuchsian), which we
denote $\Ghat$ from now on, does exist and is unique by the theorem
of Klein--Poincar\'e. The groups \Ghat\ and $\Gtilde_z$ are
constructed as  shown in
sect.\,6. If we have no some special symmetries,
the group $\Gtilde_z$
coincides with $\Gtilde$. The second generator of the field
$\boldsymbol{F}$ on Klein's curve $y=\varphi(\tau)$ in
form \reff{*} becomes a single-valued function by the construction.
Therefore $\G_x=\Ghat$. The associated equations
\reff{Q12}, including equations on initial generators
$(X,Y)$, are computed involving the {\bf Lemma}. The monodromy
$\G_y$ is Fuchsian due to the last property of the function
$\varphi(\tau)$. We have also an obvious equality $\big|\Ghat:\G_y
\big|=4$.

Parabolicity of the vertexes $\alpha$ in \reff{general} and
\reff{even} leads to  further advantage by comparison with
Whittaker's case. Namely, we do not need a birational
transformation  because arbitrary multi-valued algebraic function
$y(x)$ with ramifications only at  the points $\alpha_k^{}$ will be a
single-valued one $\varphi(\tau)=y\big(\chi(\tau)\big)$.

All the described above is a direct
generalization of Riemann--Klein's \cite{klein1, klein3}
considerations of the function $k^2(\tau)$ on the thrice punctured
sphere\footnote{One of simplest ways to uniformize Klein's curve
is to do this by a subgroup of index 7 in $\GAMMA(2)$:
$\big\{X=\sqrt[\leftroot{-1}\uproot{1}7]{k^4k'{}^2}(\tau),\;
Y=k^2(\tau)\big\}$.}.
See for example last
sentences and a figure on p.\,93 in Wirtinger's supplement (1902)
to Riemann's lectures on a hypergeometric series \cite{riemann}
or examples \reff{kk}.

Generally,  number of the accessory parameters in the
associated equations (\ref{general}--\ref{whittaker}),
\reff{even}, is greater than in an original
equation, but all the groups are free and their properties are
simpler. Summarizing, in  preceding notation, we have
finally

\noindent {\bf Theorem}. {\em Let $\X$ be an arbitrary algebraic
curve $F(X,Y)=0$ of genus $g$. Then

\begin{itemize}
\item There exists uniformization of $\X$ by a free group $\G$
being a subgroup  of a $1$-st kind free Fuchsian group $\Ghat$ of
genus zero and index $\big|\Ghat\!:\!\G\big| \geqslant
\big[\frac{2g+8}{3}\big];$

\item Associated hyperelliptic curves $\Xtilde$ and corresponding
Fuchsian equations \reff{general} have the monodromy
$\Gtilde$ of genus $\widetilde g \gtrless g$ and generate all the
representations for  the accessory parameters of both the curves
$\Xtilde$ and $\X;$

\item The group $\Ghat$ is $(2\,\widetilde g+1)$-generated
monodromical group of the equation \reff{general}   with
generators \reff{ident} explicitly described in sect.\,$6$. The
curves $\Xtilde$, their groups $\Gtilde$, genus $\widetilde g$,
and ranks of \G\ and \Gtilde\ are explicitly computed.
\end{itemize}
} \noindent {\em Proof\/}. Associated hyperelliptic curve \Xtilde\
is constructed in a complete analogy with Klein's example described
above. Taking homogeneous representation for \X\ in form
$H(u_1^{},u_2^{},u_3{})=0$ we make an arbitrary linear
transformation $u\mapsto u'\!: \;u'_j=A_{jk}^{}\,u_k^{}$. There are
only finitely many flex-points. Hence, in order to get non-singular
form of the curve \X\ with simple branch points, we have to avoid
only finitely many tangent lines. Having this done, we return to
inhomogeneous form of this new representation $\Phi(x,y)=0$ for \X\
and compute its discriminant polynomial $D(x)$. The curve \Xtilde\
is $z^2=D(x)$. The groups \Ghat, \Gtilde, genus $\widetilde g$, rank
of \Gtilde, and Fuchsian equation on the function $x=\chi(\tau)$ are
constructed  as it was done in the sect.\,6.2. All these groups are
free of defining relations and the function $\chi(\tau)$, by
construction, is one of generator of the function field on the
origin curve \X. Index of \G\ in \Ghat\ depends on the discriminant
transformation and number of sheets of \X. Lower bound
$\big[\frac{2g+8}{3}\big]$ follows from \cite{mumford}. For some
particular curves this bound can be probably decreased  to the
number $\big[\frac{g+3}{2}\big]$ being a lowest order of meromorphic
function on a curve of genus $g$. Since groups $\G$ and $\G_y$ are
free and hopfian \cite{LS}, minimal numbers of their generators are
established by pure combinatorial methods \cite{magnus} or by direct
analysis of singular points in arising Fuchsian equations. \hfill
\rule{0.5em}{0.5em}

The theorem with minor modifications is restated for Whittaker's
equation \reff{whittaker} and its monodromy.

The group \Ghat\ is very wide one
but the fact that one function uniformizes
many algebraic curves is justifiable to call the function (and its
automorphism)
$$
x=\chi(\tau),\qquad \Ghat\equiv\mbox{{\bf Aut}}\big(\chi(\tau)\big)
$$
satisfying the equations \reff{general}, \reff{even},
\reff{whittaker}, the {\em universal uniformizing function\/}
({\em universal uniformizing group $\Ghat$\/}) {\em of genus zero}
in the sense that the group $\G$ uniformizing  arbitrary algebraic
curve may be chosen as a subgroup of $\Ghat$.

\section{Consequences of the theorem}
\subsection{Modular uniformization}
If  groups $\G_x$ and $\G_y$ are of lower genera
and have the same Fuchsian {\sc ode}, i.\,e.
$$
\varphi(\tau)=\chi\Big(\mfrac{a\,\tau+b}{c\,\tau+d} \Big)
\quad \Leftrightarrow\quad\G_y =
\big(\begin{smallmatrix}a&b_{}\\c&d\end{smallmatrix}\big)^{\sss\!-1}\G_x
\big(\begin{smallmatrix}a&b_{}\\c&d\end{smallmatrix}\big),
$$
then it is natural  to speak of {\em modular equation\/} ({\em
modular uniformization}) belonging to the function $\chi(\tau)$
(group $\G_x$). From technical viewpoint such a kind of
uniformization is probably simplest  because we have
one group of a lower genus.
It is of interest to inquire whether zero  genus
modular uniformization
exists for a given curve? If the answer is negative then
what lowest admissible  genus does
exist?  Non-modularity is of course obvious.
There are lot of algebraic functions (curves) with different
ramifications at $\alpha$-points of a given universal uniformizing
function $\chi$.

If the curve has the form $F(x,y)=0$,  then this problem reduces to
searching for  a change $x \mapsto z$ such that, due to {\bf Lemma},
the following {\sc ode}
\begin{equation}\label{xyz}
[x,z]+z_x{}^{\!\!\!2}\,Q_1^{}(z,y)=Q_1^{}(x,y)
\end{equation}
has an algebraic solution $F(x,z)=0$. Such a possibility  depends on the
function $Q_1^{}$, i.\,e. on the first generator $x$.
We note that zero genus monodromy, in general,
does not necessarily entail the rational function
$Q(x)$ (see counterexample \reff{hermit}),
but an appropriate generator always exists because
meromorphic functions on \R\ of zero genus are generated
by one element (Hauptmodulus).
Therefore, in the case under consideration,
the equation \reff{xyz} may be taken in the form
\begin{equation}\label{xz}
[x,z]+z_x{}^{\!\!\!2}\,Q_1^{}(z)=Q_1^{}(x)
\end{equation}
and we are interested in its algebraic solutions.

In some
particular examples this equation
was already written down by Jacobi. Keeping
original notation we reproduce his equation \cite[p.\,132, equation (7.)]{jacobi}:
\begin{equation}\label{kl}
3\!\left(\frac{d^2\!\lambda}{dk^2}\right)^{\!\!2}-2\,\frac{d\lambda}{dk}
\!\cdot\!\frac{d^3\!\lambda}{dk^3}+
\left(\frac{d\lambda}{dk}\right)^{\!\!2} \left\{
\left[\frac{1+k^2}{k-k^3}\right]^2
-\left[\frac{1+\lambda^2}{\lambda-\lambda^3}\right]^2
\!\!\left(\frac{d\lambda}{dk}\right)^{\!\!2} \right\}=0\,.
\end{equation}
This is nothing but the equation \reff{xz} with Heun's function
$Q_1^{}$ determined by \reff{heun}. One of its algebraic solutions
is a consequence of the relation between moduli $k=k(\tau)$ and
$\lambda=k(3\tau)$ following from the relation
between $\sqrt[4]{k}(\tau)$ and $\sqrt[4]{k}(3\tau)$,
i.\,e. from the equation \reff{modular4}.
Jacobi writes this algebraic solution
in a form of polynomial of eight degree
\cite[p.\,122, equation (1.)]{jacobi},
but this polynomial is factorizable and the actual solution is
given by a simpler relation of genus one with invariant
$J=\frac{13^3}{972}$:
\begin{equation}\label{kl3}
(k-\lambda)^4=2^4 (k^3-k)(\lambda^3-\lambda)\,,
\qquad
\left\{
k=\frac{\vartheta_2^2(\tau)}{\vartheta_3^2(\tau)},\;
\lambda=\frac{\vartheta_2^2(3\tau)}{\vartheta_3^2(3\tau)}
\right\}.
\end{equation}
Other algebraic solutions of the equation \reff{kl} correspond to
higher level  modular equations. Jacobi  considered
the next nontrivial example \reff{modular8} for level 5 and gave such a
solution in several  forms
\cite[pp.\,122--123]{jacobi}. Doing again all the additional
simplifications we obtain the correct answer (genus $g=3$):
\begin{equation}\label{5}
(k-\lambda)^6=2^6 (k^3-k)(\lambda^3-\lambda)
\big(4(k\,\lambda+1)^2-(\lambda-k)^2\big)\,.
\end{equation}

We note that  the modular uniformization, in such a view, may be
naturally treated as  a transition from the ``transcendental''
(through the parameter $\tau$) integrability of Fuchsian equations
\reff{dex}  to the algebraic one for the equation \reff{xyz}. This
transition is possible  if we have correct $A$-parameters. Indeed,
since the parameter $\tau$ is a common quantity for all the Fuchsian
groups/subgroups, the general integral of the equation \reff{xyz}
has the following bilinear form
\begin{equation}\label{hyper-xz}
a\,\Psi_1(x)\Psi_1(z)+b\,\Psi_1(x)\Psi_2(z)+
c\,\Psi_2(x)\Psi_1(z)+d\,\Psi_2(x)\Psi_2(z)=0\,.
\end{equation}
Although this integral has many algebraic solutions,
in that form, it is transcendental but this transcendence
is of different kind. It is realized not through the
parametrizing functions
but through the relations between solutions of Fuchsian equations.
In particular, in the case of hypergeometric equation for the group
$\GAMMA(2)$, the formula
\reff{hyper-xz} provides generalization of classical Jacobi's
results for integer levels. We illustrate these arguments
in example of the level 3. As a consequence of the formulas
(\ref{k2}--\ref{psi}), we get
two explicit forms of the relation \reff{hyper-xz}. One of them is
\begin{equation}\label{KL}
3\,K'(k)\,K(\lambda)=K'(\lambda)\,K(k)\,,
\end{equation}
where $k$ and $\lambda$ are connected by the elliptic curve \reff{kl3}
with indicated parametrization.
Such relations, as relations between complete elliptic integrals
corresponding to different elliptic moduli, are well-known.
But here we emphasize that
quantities $k,\lambda$ are {\em not treated as elliptic moduli\/}, but
as arbitrary quantities sitting on algebraic curve of higher genus
with a given parametrization.
The curve \reff{5} is a good example on this subject.
Designating the hypergeometric function
${}_2F_1\!\!\left(\frac 12, \frac 12;1\big| x\right)$ as ${}_2F_1(x)$
we get the second form of the  identity \reff{KL}
$$
3\!\cdot\!{}_2F_{1\!}(1-\varkappa)\,{}_2F_{1\!}(\mu)=
{}_2F_{1\!}(1-\mu)\,{}_2F_{1\!}(\varkappa)\,,
$$
where variables $(\varkappa,\mu)$ belong to the reducible curve
of genus zero (not one as \reff{kl3})
$$
(\varkappa-\mu)^4=2^7(\varkappa^2-\varkappa)
(\mu^2-\mu)(2\,\varkappa\,\mu-\varkappa-\mu+2)\,.
$$
This curve, in a form of relation connecting variables $\varkappa=k^2$
and $\mu=\lambda^2$, can also be found in \cite[p.\,122]{jacobi}.
It is above mentioned polynomial of eight degree in $(k,\lambda)$.

The equations (\ref{hyper-xz}--\ref{KL})
are invariant with respect to  change of the parameter $\tau$
(as it has been eliminated) but, in general,  are not invariant
under automorphisms of a curve to which they belong.
For example, the change of parametrization using the automorphism
$(\varkappa,\mu)\leftrightarrow(\mu,\varkappa)$
is not allowable.

To summarize, we would like to remark that
once the $\Psi$-function is known, the equations
\reff{xyz} and \reff{hyper-xz} may be thought of as
{\em generating equations\/}
or the {\em universal transcendental
representation
for all curves uniformized by a zero (or higher) genus
modular uniformization
with the monodromy $\G_x$ for this function $\Psi(x)$\/}.

\subsection{Hierarchies and embedding of curves}
Depending on preferences, an initial curve can be chosen either as a
curve with relatively few branch-points, or as a cover of fairly low
degree. For the 3-branch-point covers (independently of genus and
orders of ramifications) the $A$-problem can be considered as solved
since the existence of such  representations immediately leads to a
hypergeometric uniformization by a subgroup of $\GAMMA(2)$ with
automatical producing the $A$-parameters. Four points lead to Heun's
equations and, if the Heun equation coincides with equation
\reff{heun}, the problem may be considered as completely solved
through the $\vartheta(\tau)$-constants. The hyperelliptic curves
are completely described by a universal way. In this case, all the
objects, except $A$-parameters, are explicitly constructed and the
curves themselves form a {\em hyperelliptic tower\/} of
non-isomorphic curves of descending genera $\widetilde g,\,
\widetilde g-1,\, \ldots$ up to non-isomorphic elliptic ones. All
this tower  is uniformized by  one function $\chi(\tau)$
corresponding to the genus $\widetilde g$. For example in the
simplest case of  $\widetilde g=2$ we arrive at not
$\wp$-parametrizations of  elliptic curves. Here is one example
\footnote{Comparing this $\vartheta$- and usual
$\wp$-parametrization of the lemniscate one gets new representation
for Gaussian lemniscate constant:
$\omega=\sqrt[4]{2}\,\frac{\pi}{2}\,
\vartheta_4^2(2i)=1.85407467862567819586995\ldots,\;\;
\omega'=i\,\omega$.}:

\begin{equation}\label{ell}
y^2=x^3+x\;, \qquad\qquad
\left\{
x=\frac{\vartheta_3(4\tau)}
{\vartheta_2(4\tau)},\quad
y=\frac{1}{\sqrt{8}}\,\frac{\vartheta_2^2(\tau)}{\vartheta_2^2(4\tau)}
\right\}\;.
\end{equation}
Indeed, this elliptic curve is isomorphic to the lemniscate
$y^2=x^4-1$ which, in turn, is embedded into Burnside's curve
$y^2=x(x^4-1)$. It is not difficult to see that Burnside's function
\reff{prop2} is, in fact, equal to the $x$-function in \reff{ell} up
to some linear  transformation of the global parameter $\tau$.
Another example is related to the function \reff{kkk}. It generates
completely and explicitly uniformizable tower for the curve
$y^2=z^8-1= (z^4-1)(z^4+1)=\ldots$. We also note that all the
elements of the tower fit into the {\bf Theorem} independently of
each other and generate their own hierarchies of curves and their
extensions (see below). These points are, in fact, uniformizing
realization of the classical radical functions \cite{baker} and very
closed to Poincar\'e's notion of subordination of types of Fuchsian
equations \cite[\S\,VI]{P84}.

\subsection{Mutual determination of moduli}
Since there are infinitely many curves \Xtilde\  corresponding to a
non-hyperelliptic equation \X\ and conversely, this phenomena
implies certain structures in  the space of moduli
$\boldsymbol{\mathcal M}_{g,n}$ of algebraic curves.

\noindent
{\bf Corollary 2.} {\em $A$-parameters and therefore moduli or
period-matrix of the algebraic curve
\X\ are determined by the hyperelliptic moduli of the associated curves
\Xtilde\/}.

If we  consider the $A$-parameters as moduli of the curve \X\ (see
survey of Bers \cite{bers}) then $3\,g-3$ non-hyperelliptic moduli
are computed through the $2\,\widetilde g-1$ hyperelliptic ones
({\bf Lemma}). If one deals with a transcendental representation for
moduli through $\Pi$-matrices of holomorphic integrals we should
involve, besides the hyperelliptic function $\chi(\tau)$, the second
unformizing function $\varphi(\tau)$. Hence non-hyperelliptic
$\Pi$-matrices inherit many properties of hyperelliptic ones but we
should take into account non-hyperelliptic structure of the curve
\X. Classification of such structures, i.\,e. equations of curves of
higher genera, is a highly nontrivial problem and far from being
completely solved \cite{mumford}. In Fuchsian formulation the
problem reduces, due to {\bf Lemma}, to choice of {\em the canonical
representative for the moduli\/}. The canonical representative for
all the 3-branch-point covers is unique since these curves have no
parameters. It is Legendre's equation \reff{k2} with $\GAMMA(2)$.
Four points correspond to the one $A$-parameter family of elliptic
curves uniformized by Heun's equation \reff{general} under $g=1$.
Description of moduli for such a representation of the elliptic
family is still unknown apart from one isolated point --- the
lemniscate. The solution of this problem would lead us to explicit
description of some 1-dimensional subspaces in $\mathcal{M}_{g,4}$.
Proposed zero-genus uniformization is, perhaps, a good candidate for
the problem of canonicity, even though the curve \Xtilde\ is not
unique. This is because Fuchsian equations and automorphisms $\G_z$,
under fixed genus $g$, have much more complex structure than
automorphisms of zero genera.

Some extra examples below exhibit mutual determination of the moduli
between different kind  curves. To all appearances, hierarchies of
cyclic covers $y^m=(x-\alpha_1^{})\cdots(x-\alpha_n^{})$ can  be
described effectively likewise hyperelliptic cases. Conversely,
having the hyperelliptic uniformizations we can get  their
extensions without discriminant transformation. The following two
examples suggest themselves from formulas \reff{prop2}. The first
one is
$$
z^3=x^5-x,\qquad\quad
\left\{x=\frac{\vartheta_4^{}(\tau)}{\vartheta_3^{}(\tau)}, \quad
z=-\sqrt[3]{16}\,\frac{\eta^2(2\tau)}{\vartheta_3^2(\tau)}\right\}.
$$
The function $w=\sqrt{z}(\tau)$ is  a single-valued one. We thus
obtain the second example
$$
w^6=x^5-x,\qquad\quad
\left\{x=\frac{\vartheta_4^{}(\tau)}{\vartheta_3^{}(\tau)}, \quad
w=i\,\sqrt[6]{16}\,\frac{\eta(2\tau)}{\vartheta_3^{}(\tau)}
\right\}.
$$
As in the  example \reff{ell}, the cyclic extensions are also
generated by any elements of the tower. For example, elliptic curve
\reff{ell}, belonging to this tower, yields a genus 3 cover of the
form $y^4=x^3+x$. We omit explicit formulas for its parametrization
and Fuchsian equations  with their ${}_2F_1$-solutions. Other cyclic
extensions, say $z^m=x^5-x$  or $z^m=x^3+x$, are described
completely by proposed methods except that the explicit
representation for the second generator $z(\tau)$ is unknown in
terms of known $\vartheta$-constants. In the similar manner one
obtains  {\em non\/}-modular uniformizations arising after comparing
Burnside's curve with  all the modular equations for the quantities
$\sqrt[\leftroot{1}\uproot{1}4]{k}(n\tau)$.

We also note that the list of examples can be enlarged with the help
of suitable analysis of earlier studied cases of zero genera. Most
considerable of them are recent results of Harnad \& McKay on
classification of hauptmoduli and their groups commensurable with
$\GAMMA(1)$. See tables and $\vartheta,\eta$-formulas in
\cite{mckay2,harnad,mckay,birch}, one example of computed monodromy
in \cite{hempel}, and some information on non-congruence subgroups
of $\GAMMA(1)$ in \cite{atkin}.

The associated curves may gain  nontrivial parameters coming from
the discriminant transformation even though the origin curve had no
the formers. Besides associated curves \Xtilde, we also get  curves
$\Phi(z,y)=0$  (uniformizable) obtainable by  elimination the
universal function $x$ from the origin curve $F(x,y)=0$ and the
associated one $z^2=D(x)$. The question of exact correlation between
these families of curves \Xtilde, \X, their moduli, and the
discriminant transformation \reff{*}$\mapsto$\reff{hyperKlein}
requires detailed investigation.  Another important question is how
to describe explicitly relations between  monodromies of the
generators $\G_{x,y}$ and monodromy $\G_z$ of  given a function $z$,
and how does a tessellation of \R\ corresponding to $\G_z$ look?
Answer to these questions would give answer to the question about
meromorphic functions having the same monodromy as the group \G\
uniformizing the curve itself.

\subsection{Mixed uniformizations}
If  Fuchsian equation contains both
parabolic and Whittaker's singularities, then
this also leads to the universal function, but all the groups
and polygons will have more complex structure.
One example directly follows from  \cite[\S\,4]{br}.
In notations of this work we have the following Fuchsian equation
$$
\Lambda_{\lambda\lambda}^{}= \ds -\frac14\,
\frac{\lambda^6+4\,\lambda^5+16\,\lambda^4-56\,\lambda^3+68\,\lambda^2
-48\,\lambda+16}{\lambda^2\,(\lambda-1)^2\,(\lambda^2+4\,\lambda-4)^2}^
{\ds\mathstrut}_{}\, \Lambda
$$
with three parabolic singularities $\lambda_{1,2,3}^{}=\{0,1,\infty \}$
and two ``Whittaker's'' ones $\lambda_{4,5}^{}=\{-2\pm 2\,\sqrt{2} \}$.
Every algebraic curve having  arbitrary ramification at points
$\lambda_{1,2,3}^{}$ and simple branch points at places $\lambda_{4,5}^{}$
will be uniformized by this function.
For instance, the equations
$\mu^{2g}=\lambda(\lambda-1)(\lambda^2+4\,\lambda-4)^g$
exemplify one family of such curves of genus $g$.
Universal function $\lambda(\tau)$ for all this family
reads as
$$
\lambda=\frac{\vartheta_3^2(2\tau)}
{\vartheta_4(\tau)\,\vartheta_3(4\tau)}\;.
$$
Proof is given by  the straightforward substitution any
$\vartheta$-formula for Burnside's  $\chi(\tau)$ into the formula
for $\lambda$ \cite[formula (23)]{br}, supplemented with the
subsequent simplification due to identities \reff{tt} and linear
transformation of $\tau$. Note that under $g=1$ we arrive at an
elliptic torus which is not isomorphic to the tori  \reff{ell},
\reff{kl3}, or to the torus which is covered by Burnside's curve and
described in \cite[\S\,4]{br}.

We conclude this section with  the remark that
the uniformizing equations of Whittaker, Weber and \reff{general} are
included in the general set of the same Fuchsian differential
equations (up to the multipliers), except that the accessory
polynomials will have different coefficients.

\section{Some applications and bibliographical remarks}

\noindent In order to exhibit a capacity for work of the proposed
approaches, we display here a few self-suggested applications
following from the uniformization machinery expounded in the
preceding sections. All the formulas are verified/derived by direct
substitutions and differentiating with use of {\bf Lemma}, {\sc
Appendix}, and the system \reff{fundamental}. We omit details but
accompany the material by brief literature comments. All the results
below, especially the first explicit examples of Abelian integrals
(sect.\,10.5) and their $\Theta$-functional generalizations, will be
written up later.

Schwarz's equation \reff{general} or, more precisely, its linear
Fuchsian relative, was a subject of many deep investigations
beginning with Riemann's lectures of 1858/59 \cite{riemann}. Various
single-valued $q$-series for moduli $k,\,k'(\tau)$ and their
roots/logarithms arose already in works by Jacobi \cite[{\bf I}:
pp.\,159--164]{jacobi} and Riemann \cite[pp.\,455--465]{riemann}.
See also Dedekind's comments to \cite{riemann} on pp.\,466--478. The
Fuchsian equation with $n$ parabolic singularities was a start point
of underlying ideas of Poincar\'e. Automorphic $\Theta$-series, four
parabolic singularities and Abelian integrals as functions of $\tau$
arose in \cite{poincare1} (see also dissertations
\cite{smirnov,ger}). General Fuchsian equation of the type
\reff{general}, real $A$-parameters, free monodromy $\Ghat$ (our
group $\Ghat$   differs from Poincar\'e's), and the functional
analog of Legendre's modulus arose in \cite{poincare3, poincare4}.

\subsection{Classical modular equations}
First modular equations (\ref{modular4}--\ref{modular8}) arise on
first pages of {\em Fundamenta Nova\/} \cite[\S\S\,13--15]{jacobi}
{\em before\/} the celebrated ``elliptic part'' of Jacobi's theory.
Explicit parametrization \reff{8}, in a form  of product, appeared
in a letter from Jacobi to Legendre in 1828 \cite[{\bf I}:
pp.\,409--416]{jacobi}. There seems little doubt that he considered
these and other products  as single-valued objects.

Systematic investigations of modular equations as algebraic curves
of a given genus, were started by Klein \cite[pp.\,256--273]{klein3}
and his pupils Dyck \cite[Jacobi's curve \reff{modular4}]{dyck} and
Fricke \cite{468,466,467}. Jacobi's modular equations of lower
levels lead  to relations of higher genera  whilst equations on the
invariant $J$ have genus zero  up to level 11. The explanation is
that $\GAMMA(2)$ is a free group and $\GAMMA(1)$ is not.

\subsection{Polygons and analytical representations}
Compact symmetrical form of Klein's curve, birational
transformations, the known fourteengon, and parametrization in terms
of $q$-series were written up in the famous paper of Klein
\cite{klein1} (see also \cite{haskell}). The explicit representation
all of these objects in the language of $\vartheta$-constants,
division by seven, and other results were obtained slightly later by
Klein himself \cite{klein2} (the fact rediscovered recently in
\cite{far}). It is noteworthy, however, that  Klein's
$\mathbb{H}/\GAMMA(7)$-parametrization by $\vartheta$-constant does
not correspond to his celebrated 14-gon built by $2\!\cdot\! 168$
triangles of the shape
$\big(\frac{\pi}{2},\frac{\pi}{3},\frac{\pi}{7}\big)$. Klein does
not elucidate this nontrivial point, i.\,e. relation between the
different uniformizations\footnote{Schwarz's (2,3,7)-function does
the job but $\vartheta$-efficacy is lost:  ${}_2F_1$-series
converges slowly.} and restricts himself only to comments on
transition between corresponding $\tau$-planes (see table in the
dissertation by Klein's student Haskell \cite{haskell} on p.\,8).

It is remarkable that the relation between this 14-gon (Hauptfigur
in \cite{klein1}) and its relative in $\omega$-plane ($\omega$ is
Klein's notation for $\tau$) is the same as between Whittaker's
polygon and our parabolic one. One side corresponds to two other
ones with a common parabolic point.

Somewhat surprising facet, but it should be confessed, is the fact
that the analytical description in the language of Fuchsian {\sc
ode}'s has not got to monumental treatises
\cite{483,467,468,539,496,lehner,farkas2000,551} and later works on
uniformization. This especially concerns a reducibility of the
monodromies to ones of lower genera (see however
\cite[\S\,XVI]{P84}). Whittaker was the first to joint explicitly
polygons, Fuchsian equations \cite{418}, and reductions of groups
\cite{420}. Even though Poincar\'e was interested in generalization
of Legendre's modulus \cite{poincare4}, applicability of the
equations (\ref{general}--\ref{whittaker}) to the cases of higher
genera is not clear from correspondence between him and Klein
\cite[pp.\,587--\ldots]{klein3} or from their later works (arguments
in last paragraph of \cite[\S\,XVII]{P84} are very unclear). This
is, probably, due to purely zero genus Klein's treatments of maps of
half-planes onto triangles (see pages 592, 597, and 600 in
\cite{klein3}).

\subsection{Uniformization and solution of algebraic equations}
One direct application  of modular equations is a construction of
resolvents and their solutions \cite{weber}. Relationship of
uniformization with this topic was, in implicit form, mentioned in
the literature, but has never been realized explicitly. Some known
examples do not use ideas of uniformization (see for example
\cite{knopp}). Hermitian resolvent (the quintic) is not an exception
in spite of much similarity of approaches. The base idea in such a
view is the fact that  roots of the equations are equivalent points
belonging to different sheets of \R.

As well-known, the general quintic depends on one essential
parameter $a$ and, after suitable Tschirnhaus's transformations
(E.\,S.\,Bring 1786), can be brought into the canonical form
$x^5-x-a=0$. Not invoking the resolvent view, we make  use of
Burnside's parametrization and get the answer\footnote{It is
curiously that the celebrated Hermitian solution of the quintic
(1858)

$$
x^5-x-\mfrac{2}{\sqrt[4]{5^5}}
\mfrac{1+\varphi^8(\omega)}{\varphi^2(\omega)\,\psi^4(\omega)}=0,
\qquad x_k^{}=\!\!\sqrt[-4]{2^45^3}\,
\mfrac{\Phi(\omega+16k)}{\varphi(\omega) \,\psi^4(\omega)} \qquad
(k=0\ldots 4)
$$
contains an erroneous second sign ``$-$'' independently of
definitions of
$\varphi(\omega)=\frac{\vartheta_4^{}(2\omega)}{\vartheta_3^{}(\omega)}$
and, to the best of our knowledge, this misprint has crept into all
subsequent reprints of the solution. Clearly, uniformization catches
automatically all the similar things coming from  a residue symbol
or $\sqrt[8]{1}$'s. It also explains further Hermitian emphasis
``{\em le produit de deux fonctions $\varphi(\omega)$ et
$\psi(\omega)$ est le cube d'une nouvelle fonction \'egalement bien
d\'etermin\'ee\/}'' ({\sc C.\,Hermite:} {\em \OE uvres\/} {\bf II}
(1908), p.\,28).}

$$
x^5-x+2^4\,\frac{\eta^6(2\tau)} {\vartheta_3^6(\tau)}=0,
\qquad\qquad x_{1\ldots 5}^{}= \frac{\vartheta_4^{}(\tau_{\sss
1\ldots 5})} {\vartheta_3^{}(\tau_{\sss 1\ldots 5})},
$$
where $\tau_{1\ldots 5}^{}$ are images of the point $\tau$ in a
tessellation of the  18-gon ({\sc Fig.\,2}) by five $y$-sheets built
on $\tau$-pre-images of the branch points $y(\tau_k^{})=\big\{0,
\;2\!\sqrt[-8]{5^5},\;\infty \big\}$. The sheets $x_{1\ldots 5}^{}$
have to be interchanged by the group $\G_{y^2}$ of the function
$$
y^2=-2^4 \,\frac{\eta^6(2\tau)}{\vartheta_3^6(\tau)}\qquad (\equiv
z).
$$
Its monodromy (genus zero) is governed by the following Fuchsian
equation:
\begin{equation}\label{hermit}
\widetilde\Psi_{zz}=-\frac14\,
\frac{z^8-5^{-3}48\,x\,z^7+12^25^{-4}\,x^2z^6+\cdots+2^{16}5^{-10}}
{z^2(z^4-2^85^{-5})^2}\,\widetilde\Psi.
\end{equation}
Mentioned tessellation\footnote{The problem of the similar
tessellations has not been solved even in the elliptic case. This
also explains why we marked the sign $\pm$ at point 5 in the table
of sect.\,1.2.} of \R\ is of interest in its own right and can be
obtained by inversion procedure described in sect.\,5.1 or making
use of Hermitian formulas for the quintic (see preceding footnote).
Sought for tessellation is then obtained with the help of 5-fold
Hermitian transformations the circular arcs bordering the 18-gon \P\
in {\sc Fig.\,2}. Ultimate formulas for this solution of the quintic
will be presented elsewhere.

The similar approach can be applied to different forms of quintics
obtainable from other parameterizations. For example from our curve
\reff{new} or Weber's example in \cite[{\bf III}, p.\,256]{weber},
etc. See also \cite[pp.\,366--376]{lehner},
\cite[pp.\,119--123]{knopp} and \cite{birch} for $\GAMMA_0(50)$.
Efficacy of formulas depends on a geometry of polygons.

\subsection{Jacobi--Chazy- and Picard--Fuchs-like equations}
First differential relation satisfied by some modular object was the
equation \reff{kl} between moduli $k(\tau)$ and $k(n\tau)$. Later,
Jacobi derived an important third order {\sc ode} for the series
$$
y=1\pm 2\,q+2\,q^4 \pm 2\,q^9 +\cdots,
$$
i.\,e. for $y=\vartheta_{3,4}^{}(\tau)$-constants (Jacobi claims
this result as {\em Theorema}  \cite{jacobiCrelle}):
$$
\big(y^2\,y'''-15\,y\,y'\,y''+30\,y'{}^3\big)^2
+32\,\big(y\,y''-3\,y'{}^2
\big)^3=-\pi^2\,y^{10}\big(y\,y''-3\,y'{}^2 \big)^2\;.
$$
Analog of this
equation for modular discriminant $\Delta(\tau)$ was
rediscovered many times later (van der Pol (1951),
Res\-ni\-koff (the 1950--60's), Ablowitz et all (the 1990's),
\cite{494, takh2}). In the modern language, these equations are
obviously none other than equations on
modular forms. They have an
independent interest and theory, but are direct consequences of
Fuchsian equations (\ref{dex}--\ref{Q12}).
General recipe for obtaining such formulas is as follows.

The function $\Psi^2=x_\tau$ is an automorphic form of weight 2 and,
hence,  satisfies some differential equation of third order
(Hurwitz). Rewriting the equation \reff{dex} and its $x$-derivative
in variable $\tau$, we get the following two equations:
$$
2\,\Psi\,\Psi_{\tau\tau}-4\,\Psi_\tau{}^{\!\!\!2}=Q(x)\,\Psi^6,\qquad
2\,\Psi^2\Psi_{\tau\tau\tau}-18\,\Psi\,\Psi_\tau\Psi_{\tau\tau}
+24\,\Psi_\tau{}^{\!\!\!3}=Q_x(x)\,\Psi^9.
$$
The sought for formulas are obtained by elimination  $x$
from these equalities.  In this context we would like to emphasize
that if some uniformization is described by Jacobi's constants
then all relevant differential relations/identities
are nothing but consequences of one
system of {\em closed differential equations on the
$\vartheta,\eta_{\mbox{\tiny\sc w}}^{}$-constants\/}
\begin{equation}\label{fundamental}
\begin{array}{ll}
\ds
\frac{1}{\vartheta_{2\!}^{}}\,\frac{d\vartheta_{2\!}^{}}{d\tau}
_{\ds\mathstrut}= \mfrac{i}{\pi}\,\eta_{\mbox{\tiny\sc
w}}^{}+\mfrac{\pi i}{12}\,
\big(\vartheta_{3\!}^4+\vartheta_{4\!}^4 \big),&
\ds\;\;
\frac{1}{\vartheta_{4\!}^{}}\,\frac{d\vartheta_{4\!}^{}}{d\tau}
^{\ds\mathstrut}= \mfrac{i}{\pi}\,\eta_{\mbox{\tiny\sc
w}}^{}-\mfrac{\pi i}{12}\,
\big(\vartheta_{2\!}^4+\vartheta_{3\!}^4 \big),
\\
\ds
\frac{1}{\vartheta_{3\!}^{}}\,\frac{d\vartheta_{3\!}^{}}{d\tau}
_{\ds\mathstrut}^{\ds\mathstrut}=
\mfrac{i}{\pi}\,\eta_{\mbox{\tiny\sc w}}^{}+\mfrac{\pi i}{12}\,
\big(\vartheta_{2\!}^4-\vartheta_{4\!}^4 \big),&
\ds\quad\;\;\,
\frac{d\eta_{\mbox{\tiny\sc
w}}^{}}{d\tau}=\mfrac{i}{\pi}\Big\{2\,\eta_{\mbox{\tiny\sc
w}}^2-\mfrac{\pi^4}{144}
\big(\vartheta_{2\!}^8+\vartheta_{3\!}^8+\vartheta_{4\!}^8 \big)
\Big\},
\end{array}
\end{equation}
wherein $\eta_{\mbox{\tiny\sc w}}^{}$ is Weierstrass's function
$\eta_{\mbox{\tiny\sc w}}^{}=\zeta(1|1,\tau)$. This even-dimensional
version of the celebrated Halphen's system is not seem to be in the
literature. See also the handbook \cite{br2} of new formulas in
theory of elliptic/modular functions. The Halphen system itself and
its modern varieties \cite{harnad,mckay2,takh2} are consequences of
this system and algebraic identities between $\vartheta$-constants
({\sc Appendix}). If one has Dedekind's $\eta$-constant, we should
involve the equations
$\pi\,\eta_\tau^{}=i\,\eta\,\eta_{\mbox{\tiny\sc w}}^{}$ and
$$
\big\{\Lambda_{\tau\tau\tau}-12\,\Lambda_{\tau\tau}\,\Lambda_{\tau}
+16\,\Lambda_{\tau}{}^{\!\!\!\!3} \big\}^2
+32\,\big\{\Lambda_{\tau\tau}-2\,\Lambda_{\tau}{}^{\!\!\!\!2}\big\}^3=
\mfrac{4}{27}\,\pi^6\,\eta^{24}\;.\qquad (\Lambda\equiv \ln\eta)
$$

Main motivation of Poincar\'e, in his correspondence with Fuchs and
Klein in  1880/81, was not uniformization but integration of
equations of Fuchsian class in terms of new single valued functions.
See also Poincar\'e's emphasis at the end of \S\,XI in \cite{P84}.
The results of preceding sections allow us to get  such examples in
explicit form.

Under the fixed  monodromy, all the $\Psi$-functions, corresponding
to subgroups, are merely multiplied by a factor. In other words,
these are proportional to one $\tau$-representation of the {\em
universal\/} $\Psi$-function. $\GAMMA(2)$-representation for the
universal $\Psi$-function is the formula \reff{Ktheta}. Burnside's
$\GAMMA(4)$-representation is given by the  formula
\begin{equation}\label{last}
\Psi_1(\tau)=\pm 2\,i \sqrt{i\pi}\,\frac{\eta^3(2\tau)}
{\vartheta_3^{}(\tau)}\;,
\qquad \Psi_2=\tau\,\Psi_1,
\end{equation}
and the $\GAMMA(1)$-version is described by the following
equation and its solution
$$
J(J-1)\,\mfrac{d^2}{dJ^2}\eta^2(\tau)+
\mfrac16\,(7J-4)\,\mfrac{d}{dJ}\eta^2(\tau)
+\mfrac{1}{144}\,\eta^2(\tau)=0\;.
$$
Higher genera $\widetilde\Psi$-functions have the same
$\tau$-representations as \reff{last}:
$\widetilde\Psi=(A\tau+B)\Psi_1(\tau)$. In this case we have to make
a change $x\mapsto y$ and supplement it with the linear change
$\Psi(x)=\Psi\big(x(y)\big)=\widetilde\Psi(y)$ without scale factor
$\sqrt{y_x^{}}$.

Among other single-valued objects of the theory we also exhibit
explicit $\tau$-formulas for Schwarz--Klein's octahedral
ground-forms \cite{klein3,483,br} ($\vartheta\equiv
\vartheta(\tau)$)
$$
\begin{array}{rcr}
\vartheta_4^{}\vartheta_3^{}\big(\vartheta_4^4-\vartheta_3^4
\big)\!\!\!&=&\!\!\!
-16\,\eta^6(2\tau)\;,\\
\vartheta_4^8+14\,\vartheta_4^4\,\vartheta_3^4+\vartheta_3^{8\ds\mathstrut}
\!\!\!&=&\!\!\! 3\,\mfrac{2^{6\ds\mathstrut}}{\pi^4}\,g_2^{}(2\tau)\;,\\
\vartheta_4^{12}-33\,\vartheta_4^8\,\vartheta_3^4 -33\,
\vartheta_4^4\,\vartheta_3^8+\vartheta_3^{12\ds\mathstrut}\!\!\!&=&\!\!\!
-3^3\mfrac{2^{9\ds\mathstrut}}{\pi^6}\,g_3^{}(2\tau)\;\,
\end{array}
$$
and some examples of solvable linear {\sc ode}'s with modular
(and the like) coefficients:
$$
\psi_{\tau\tau}^{}+ \left\{
\mfrac{\pi^2}{288}\Big(7\vartheta_4^4+
\vartheta_3^4+\mfrac{48}{\pi^2}\,\eta_{\mbox{\tiny\sc
w}}^{}(\tau)\Big)^2- \mfrac{3}{\pi^2}\,g_2^{}(2\tau)
 \right\}\psi=0, \qquad \psi=\Psi_1(\tau),
$$
$$
\psi_{\tau\tau}^{}+\mfrac{\pi^2}{4}\,\vartheta_4^8(2\tau)\,\psi=0,\qquad
\psi=A\,\vartheta_2^{-2}+B\,\vartheta_2^{-2}
\ln\frac{\vartheta_4}{\vartheta_3}\,.
$$

\subsection{Abelian integrals}
From the uniformization viewpoint, basic Abelian integrals as
functions of $\tau$, should be considered as independent
single-valued objects. If, however, the curve can be realized as a
cover of a lower genus curve then  the integrals are written down
explicitly in terms of lower genera integrals and functions up to
Weierstrassian ones. For example, with use of formulas
(\ref{vartheta}--\ref{tt}) and results of \S\,4 in work \cite{br},
we obtain the following form for basic holomorphic integrals in
representation for Burnside's parametrization:
\begin{equation}\label{holo}
\begin{array}{l}
\ds\int\limits^{\,\,x}\!\frac{x\mp i\sqrt{i\,}}{\sqrt{x^5-x}}\,dx=
-\pi \int\limits^{\,\,\tau}\!\!\big\{\vartheta_4(\tau)\mp
i\sqrt{i\,}\,
\vartheta_3(\tau)\!\big\}\,\eta^3(2\tau)\,d\tau=\\
\phantom{\ds\int\limits^{\,\,x}\!\frac{x\mp
i\sqrt{i\,}}{\sqrt{x^5-x}}\,dx} = -\sqrt{1+i\,}\,\wp^{\mbox{-}1}\!\!
\left(\mfrac{\big(1\pm\sqrt{2}\big)\,\vartheta_3^2(2\tau)}
{2\,\vartheta_4(\tau)\,\vartheta_3(4\tau)}-\mfrac{3\pm\sqrt{2}}{6};
\,\mfrac53,\,\mfrac{\mp7}{27}\sqrt{2}\right).
\end{array}
\end{equation}
See \cite[\S\,5]{br} for {\em single-valued\/} $q$-series
representations for these holomorphic functions. In what concerns
non-holomorphic integrals, we shall exhibit only two integrals with
single poles at branch places $x=i,\,-1$ ($\tau=2$ and $+i\infty$
respectively):
$$
\begin{array}{l}
\ds I_1=\int\limits^{\,\,x}\!\!\frac{1}{(x-i)}\frac{dx}{\sqrt{x^5-x}}\,=
-(1+i)\,\frac{\pi}{32}\,
\int\limits^{\,\,\tau}\!\frac{\vartheta_4(2\tau)\,\vartheta_2^4
\big(\frac{\tau}{2}\big)} {\vartheta_4\big(\tau+\frac12\big)}\,
\,d\tau,\\
\ds
I_2=\int\limits^{\,\,x^{\ds\mathstrut}}\!\!\frac{1}{(x+1)}
\frac{dx}{\sqrt{x^5-x}}=\phantom{(1+i)}
\!-\frac{\pi}{32}\,
\int\limits^{\,\,\tau}\!\frac{\vartheta_4(2\tau)}{\vartheta_3(4\tau)}\,
\vartheta_2^4\Big(\mfrac{\tau}{2}\Big)\,d\tau.
\end{array}
$$
Explicit formulas to them are extracted from the following two
linear identities
$$
\mfrac{1\pm\sqrt{i\,}}{\sqrt{1+i\,}}\,I_1+
\mfrac{1\mp i\sqrt{i\,}}{\sqrt{1+i\,}}
\,I_2=
\zeta\Big(\alpha^{\pm}(\tau);
\mfrac53,\mfrac{\mp7}{27}\sqrt{2}\Big)+
\mfrac{3\pm2\sqrt{2}}{6}\,\alpha^{\pm}(\tau)-
\mfrac{2\pm\sqrt{2}}{2}\,\alpha^{\mp}(\tau),
$$
where $\alpha^{\pm}(\tau)$ are proportional to the two holomorphic Abelian
integrals \reff{holo}
$$
\alpha^{\pm}(\tau)=\wp^{\mbox{-}1}\!\!
\left(\mfrac{\big(1\pm\sqrt{2}\big)\,\vartheta_3^2(2\tau)}
{2\,\vartheta_4(\tau)\,\vartheta_3(4\tau)}
-\mfrac{3\pm\sqrt{2}}{6};
\,\mfrac53,\,\mfrac{\mp7}{27}\sqrt{2}\right).
$$
Explicit transition between the curve \reff{bc} and elliptic tori
$(\alpha^{\pm})$ or, which is the same,
transcendental representations (covers) for Burnside's curve
itself, have the form \cite{br}
\begin{equation}\label{trans}
F(\alpha^{\pm},x)=0:\quad
\left\{
\begin{array}{l}
\;\wp\Big(\alpha^{\pm};
\mfrac53,\mfrac{\mp7}{27}\sqrt{2}\Big)=
\mfrac{(1\pm\sqrt{2})(1-i)\,x}{(x-i)(x+1)}-\mfrac{3\pm \sqrt{2}}{6}\\
\wp'\Big(\alpha^{\pm}; \mfrac53,\mfrac{\mp7}{27}\sqrt{2}\Big)=
2\,\mfrac{\sqrt{1-i\,}}{\sqrt{2}\mp 2}\, \mfrac{\big(x\pm
i\sqrt{i\,}\big)\sqrt{x^5-x\,}}{(x-i)^2(x+1)^2} ^{\ds\mathstrut}
\end{array}
\right.\;.
\end{equation}
Corresponding holomorphic differentials are related as follows
$$
\mfrac{x\mp
i\sqrt{i\,}}{\sqrt{x^5-x}}\,dx=\sqrt{1+i\,}\,d\alpha^{\pm}.
$$
Analogous formulas are written down for many curves described
in the work. For example,
works \cite{amstein,baker} contain exhaustive information
for Jacobi's curve \reff{jacobi}.
Clearly, all the elliptic integrals on covered tori lift to
single-valued integrals on \R\  for \reff{bc}
and conversely, these are only cases of such reductions.
Integrals of exact differentials
(meromorphic functions or logarithms of them) are of
$\vartheta$-constants and expressible in terms of
$\vartheta$-constants themselves. Among other things,
such formulas can be of interest in their own right as new nontrivial
integral identities on $\vartheta$-constants.
Abelian integrals as functions of $\tau$ certainly requires the
detailed description.

\subsection{Liouville equation and metrics of Poincar\'e}
Besides mathematical interest, the equations
(\ref{general}--\ref{whittaker}) have numerous and very deep
physical motivations from modern  conformal field theories and, as
we have seen now, their universality is extended over nontrivial
situations of higher genera. In many respects, the space
$\mathcal{M}_{g,n}$ inherits the structure and moduli of  orbifolds
corresponding to zero genera monodromies.

\noindent
{\bf Corollary 3.} {\em Functional solution of the $A$-problem
obtained by Zograf and Takhtajan for the equation \reff{general}
is carried over to Riemann surfaces of higher genera\/}.

Advantages of such an approach were pointed out above:
(1) all groups are free; (2) hierarchies of curves and
mutual determination of moduli; (3) explicit examples.
No doubt that proposed computational aspects  will provide an
additional efficacy to these theories and quantum
Liouville equation intensively studied last decades in works by
L.\,Takhtajan et all \cite{tah3, tah5,tah4,takh1}.  Here are few examples.

Poincar\'e showed
that   metric of a constant  curvature $-1$ on the
universal cover $\mathbb{H}$:
$$
dS^2=\frac{|d\tau|^2}{(\Im\tau)^2},
$$
induces the conformal metrics $dS^2=e^{2U}dx\,dx^*$ on \X.
They are governed by
the complex-non-analytic but real-valued  function $U(x,x^*)$
satisfying two differential equations \cite{pl}: nonlinear {\em
non-analytical\/} {\sc pde} of Liouville and {\em analytical\/}
{\sc ode} in $x$:
\begin{equation}\label{liouville}
4\,U_{xx^*}=e^{2U}_{},\qquad
{}-U_{\mathit{xx}}+U_x{}^{\!\!\!2}=\frac12\,Q_1^{}(x,y(x)).
\end{equation}
The last equation is equivalent to the linear equation \reff{Q12}
\cite{pl}
$$
\big(e^{-
U}\big)_{\!\mathit{xx}}=\frac12\,Q_1^{}(x,y)\,\big(e^{-U}\big)
$$
and we observe that all these equations are equivalents  of each
other in the sense that the solution $U(x,x^*)$ is realized in
a form of {\em separability\/} of variables $(x,x^*)$\footnote{The Liouville
equation \reff{liouville} has also Lie point symmetry with
separable variables in its  generator $\,\widehat{\boldsymbol{\!
G}}=\tau\,\partial_x+\tau^*\partial_{x^*}-
\frac12\big(\tau_x+\tau^*_{x^*}\big)\partial_U^{}$ wherein $\tau$
is an arbitrary analytic function of $x$.} through the linear
analytical equation of Fuchs \reff{Q12}:
$$
e^{-U(x,x^*)}=A(x^*)\,\Psi_1(x)+B(x^*)\,\Psi_2(x)\,.
$$
Different normalizations for $A(x^*),\,B(x^*)$  depend only on a
model of the universal cover. Thus all the conformal metrics
$$
e^{2U(x, x^*)}=\frac{-4\,|\tau_x|^2}{(\tau- \tau^*)^2}
$$
solve the equations \reff{liouville} with an arbitrary analytic
function $\tau=\frac{\Psi_2(x)}{\Psi_1(x)}$ and, among them, the
metrics of Poincar\'e correspond to  Riemann surfaces of {\em
finite\/} genus. This requires the correct choice of the
$A$-parameters in equation \reff{Q12} on the functions
$\Psi_{1,2}(x)$ determining the ratio $\tau(x)$. The metrics have
the form
\begin{equation}\label{metrics}
\begin{array}{rl}
\boldsymbol{\Delta}/\G_x:&\quad
dS^2=\mfrac{4\,dx\,dx^*}{\big(|\Psi_1|^2-|\Psi_2|^2\big){}^2}\,,\;\,
\quad
\big\{A(x^*)=\Psi_1^*,\;\;\;\;\;B(x^*)=-\Psi_2^*\big\}\\
\mathbb{H}/\G_x:&
\quad dS^2=\mfrac{-4\,dx\,dx^{*{}^{\ds\mathstrut}}}
{\big(\Psi_2^*\Psi_1^{}-\Psi_1^*\Psi_2^{}\big){}^2}\,, \quad
\Big\{A(x^*)=\mfrac{i}{2}\,\Psi_2^*,\;\;B(x^*)=-
\mfrac{i}{2}\,\Psi_1^*\Big\},
\end{array}
\end{equation}
where choice of  $\Psi_{1,2}(x)$ should be strictly controlled by
the condition that the ratio $\tau$ belongs, independently of the
value $x$, to the unit disk $\boldsymbol{\Delta}$ (or $\mathbb{H}$).
Finding this unit circle (real axis) is  a question of normalization
of the $\Psi$.

In order to get  the metrics one needs to have solved Fuchsian {\sc
ode}'s with monodromies of higher genera. It is not sufficient to
have the $A$-parameters. Of all  examples, which suggest themselves
from preceding sections, we exhibit only one for 5-sheeted
Burnside's Riemann surface $\mathbb{H}/\G_y$ in coordinates
$(y,y^*)$. Making use of  formulas \reff{psi} and \reff{metrics} we
get
$$
dS^2=\frac{4\,\pi^2}{\Re\mathfrak
e^2\big\{K({x^*}^2){\ds
K'}(x^2)\big\}}\,\frac{dy\,dy^*}{|5y^3/x+4y|^2}\,.
$$
The complete elliptic integrals $K$ and $K'$, as above,
are expressible via hypergeometric function
$K(z)=\frac\pi2\,{}_2F_1\big(\frac12,\frac12;1\big|z^2\big)$
\cite{watson} and $x=x(y)$ is a root of the equation $x^5-x=y^2$.
So long as the example is an exact one, it is a good
exercise  to compare it with all ingredients of an approach to
uniformization through a regularized action functional for the
Liouville equation developed in \cite{tah5}.

We note that ``reducibility'' to larger groups leads to the lower
genera $\R$'s (or orbifolds) and their metrics (for example for
Whittaker's curves $\boldsymbol{\Delta}/\G_y$ \cite{420}) because
genus of the metric depends on a choice of the plane under covering.
However hypergeometric reduction is a simplest one to non-elementary
Fuchsian groups acting on $\mathbb{H}$. This means that formulas for
{\em nontrivial Poincar\'e's metrics containing ${}_2F_1$ are
simplest of possible\/}.

Metric for the universal uniformizing group $\Ghat$ in form
\reff{metrics} induces metrics for all its subgroups. If the
monodromy $\G_y$ coincides the group $\G$ of algebraic curve
\reff{1}, then the Poincar\'e metric for this Riemann surface
$\mathbb{H}/\G$, independently of kind of uniformization
\reff{general} or \reff{whittaker}, is given by the formula
\begin{equation}\label{Fxy}
dS^2=\left|\frac{F_y(x,y)}{F_x(x,y)}\right|^2
\frac{-4\,dy\,dy^{*{}^{\ds\mathstrut}}}
{\big(\Psi_2^*\Psi_1^{}-\Psi_1^*\Psi_2^{}\big){}^2}, \qquad
\frac{\Psi_2(x)}{\Psi_1(x)}\in \mathbb{H}\quad \forall\, x\in
\overline{\mathbb{C}}
\end{equation}
where $\Psi_{1,2}$ are determined through the universal uniformizing
function: $\Psi_1^2(x)=\chi_\tau^{}(\tau)$. We remark however that
due to the {\bf Lemma} and  the {\bf Theorem}, all the metrics on
\R's are nothing more than different {\em coordinate representations
of the universal metric of the zero genus \reff{metrics}, i.\,e.
metrics of orbifolds defined by equations \reff{general} or
\reff{whittaker}}. If the representation for \R\ has an algebraic
form as (1) then we have \reff{Fxy}. If it is transcendental, as in
the example \reff{trans}, we arrive at a transcendental
representation for the metric. In our main example \reff{trans},
both the tori $(\alpha^{\pm})$ are isomorphic and we take only one.
As it follows from \reff{trans} and \reff{metrics}, after some
simplification, we get the following representation for Poincar\'e's
metric in coordinates $(\alpha,\alpha^*)$ of this 2-sheeted covered
torus:
$$
dS^2=
\mfrac{\pi^2\,6^{-4}(3+\sqrt{2})\,\big|6\wp(\alpha)+3-\sqrt{2}\big|^8\,
\big|\wp'(\alpha)\big|^2\,d\alpha\,d\alpha^*}
{\big|(6\wp(\alpha)+5\sqrt{2}-3) D(\alpha)-
24i(1-\sqrt{2})(3\wp(\alpha)+\sqrt{2})\big|^2
\,\Re\mathfrak
e^2\big\{K({x^*}^2)K'(x^2)\big\}}\;,
$$
where
$$
\begin{array}{rcl}
D(\alpha)\!\!\!&=&\!\!\!
\pm\sqrt{\mbox{\footnotesize
$\ds\big(6\wp(\alpha)-3+\sqrt{2}\big)
\big(6\wp(\alpha)+21+13\sqrt{2}\big)$}}\;,
\\\\
\ds
x^2\!\!\!&=&\!\!\!
\mfrac{24i(1+\sqrt{2})(3\wp(\alpha)-\sqrt{2})-
(6\wp(\alpha)-5\sqrt{2}-3) D(\alpha)}
{(6\wp(\alpha)+3+\sqrt{2})^2}\;.
\end{array}
$$
This example is not without interest because the transcendental
representations for metrics through the coordinates on a curve of
lower genus are not always possible but only if the curve admits a
transcendental cover.

\section{Acknowledgements}
\noindent The author thanks Professor J.\;C.\;Eilbeck for much
attention to the work and useful comments and Professor E.\,Previato for
numerous and valuable discussions. The work was
supported by a Royal Society/NATO Fellowship and NSF/NATO-grant
DGE--0209549.

\section{Appendix}

\noindent Below we display notations and base properties of
functions used in the paper.
We designate the values of Jacobi's $\theta(z|\tau)$-functions at $z=0$
as $\vartheta$-constants:
$$
\mbox{\small$
\ds
\vartheta_2^{}(\tau)=2\,e^{\frac{\pi}{4}i\tau}_{\mathstrut}\,
\sum\limits_{k=0}^{\infty}\, e^{(k^2+k)\pi i\tau},\quad
\ds
\vartheta_3^{}(\tau)=1+2\sum\limits_{k=1}^{\infty}\,
e^{k^2\pi i\tau},\quad
\ds
\vartheta_4^{}(\tau)=1+2\sum\limits_{k=1}^{\infty}\,
(-1)^k e^{k^2\pi i\tau}.
$}
$$
We enumerate some identities  between $\vartheta$-constants which
were used in the work:
$$
\begin{array}{l}
\vartheta_2^2\Big(\mfrac\tau2 \Big)=
2\,\vartheta_2(\tau)\,\vartheta_3(\tau),\quad
\vartheta_3^2\Big(\mfrac\tau2 \Big)^{\ds\mathstrut}=
\vartheta_3^2(\tau)+\vartheta_2^2(\tau),\quad
\vartheta_4^2\Big(\mfrac\tau2 \Big)^{\ds\mathstrut}=
\vartheta_3^2(\tau)-\vartheta_2^2(\tau),
\end{array}
$$
\begin{equation}\label{vartheta}
\begin{array}{l}
2\,\vartheta_2^2(2\tau)=\vartheta_3^2(\tau)-\vartheta_4^2(\tau)
\qquad\quad\;\;\,\,=4\,\vartheta_2(4\tau)\,\vartheta_3(4\tau),\\
2\,\vartheta_3^2(2\tau)=
\vartheta_3^2(\tau)+\vartheta_4^2(\tau) \qquad\quad\;\;\,\,=
2\,\vartheta_4\big(\tau+\frac12 \big)\,
\vartheta_3\big(\tau+\frac12 \big)^{\ds\mathstrut},\\
2\,\vartheta_4^2(2\tau)=\vartheta_3^2\big(\tau+\frac12\big)+
\vartheta_4^2\big(\tau+\frac12\big)^{\ds\mathstrut}=
2\,\vartheta_4(\tau)\,\vartheta_3(\tau)^{\ds\mathstrut},\\
\end{array}
\end{equation}
\begin{equation}\label{tt}
\begin{array}{l}
\vartheta_4(\tau)+\vartheta_3(\tau)=
\phantom{-}2\,\vartheta_3(4\tau),\\
\vartheta_4(\tau)-\vartheta_3(\tau)=
-2\,\vartheta_2(4\tau),\phantom{\big|}^{\ds\mathstrut}_{}
\end{array}
\qquad
\begin{array}{l}
\vartheta_4(\tau)+i\,\vartheta_3(\tau)=
(1+i)\,\vartheta_3\big(\tau+\frac12\big),\\
\vartheta_4(\tau)-i\,\vartheta_3(\tau)=
(1-i)\,\vartheta_4\big(\tau+\frac12\big)^{\ds\mathstrut}.
\end{array}
\end{equation}
The identities containing  $\vartheta(\tau+\frac12)$  are new.
Dedekind's functions $\eta$ and $\eta_{1,2,3}^{}$ are defined as
follows \cite{dedekind}:
$$
\ds\eta(\tau)= e^{\frac{\pi i}{12}\tau}_{\mathstrut}
\,
\prod\limits_{k=1}^{\infty}\! \big(1-e^{2k\pi i \tau}\big),
$$
$$
\eta_1^{}(\tau)\equiv\eta(2\,\tau),\quad
\eta_2^{}(\tau)\equiv\eta\Big(\mfrac{\tau}{2}\Big),\quad
\eta_3^{}(\tau)\equiv\eta\Big(\mfrac{\tau+1}{2}\Big).
$$


\begin{thebibliography}{99}

\bibitem{acta} Acta Aplicand\ae\ Mathematic\ae\ (1994), {\bf 36}
(whole issue).

\bibitem{449}\mbox{\sc Ahlfors, L.\,V. \&  Sario, L.}
{\em Riemann Surfaces\/}. Princeton Univ.\;Press (1974).

\bibitem{aigon}\mbox{\sc Aigon, A. \& Silhol, R.}
{\em Hyperbolic hexagons and algebraic curves in genus 3\/}.
J.\,London Math.\,Soc. (2) (2002), {\bf 66}(3), 671--690.

\bibitem{amstein} \mbox{\sc Amstein, H.}
{\em Fonctions Ab\'eliennes du Gendre $3$\/}.
Bulletin de la Soci\'et\'e Vaudoise Sciences Naturelles
(1888), {\bf XXIV}(99), 1--74.

\bibitem{atkin} \mbox{\sc Atkin, A. \& Swinnerton-Dyer, H.}
{\em Modular forms on non-congruence subgroups\/}.
In: Combinatorics (Proc.\,Sympos.\,Pure Math.), {\bf XIX},
Univ. California 1968. AMS, Providence, R.I. (1971), 1--25.

\bibitem{451} \mbox{\sc Appell, P. \&  Goursat, E.}
{\em Th\'eorie des fonctions alg\'ebriques d'une variable.
{\bf II}$:$ Fonctions automorphes\/}.
(P.\;Fatou). Gauthier--Villars, Paris  (1930).

\bibitem{baker} \mbox{\sc Baker, H.\,F.}
{\em Abelian functions. Abel's theorem and the allied theory of theta
functions\/}. Cambridge Univ.\;Press (1897).

\bibitem{511} \mbox{\sc Beardon, A.\,F.} {\em The Geometry of Discrete Groups\/}.
Springer--Verlag (1983).

\bibitem{bers} \mbox{\sc Bers, L.} {\em Uniformization, moduli, and Kleinian groups\/}.
Bull.\;London Math.\;Soc. (1972), {\bf 4}(3), 257--300.

\bibitem{birch} \mbox{\sc Birch, B.\,J.}
{\em Some calculations of modular relations\/}.
In \cite{kuyk}, pp.\,175--186.

\bibitem{br} \mbox{\sc Brezhnev, Yu.\,V.}
{\em Uniformization$\,:$ on the Burnside curve
$y^2=x^5-x$\/}.\\
{\tt http://arXiv.org/math.CA/0111150}.

\bibitem{br2} \mbox{\sc Brezhnev, Yu.\,V.}
{\em On functions of Jacobi and Weierstrass (I)\/}.\\ {\tt
http://arXiv.org/math.CA/0601371}.

\bibitem{303} \mbox{\sc Burnside, W.}
{\em On a class of automorphic functions\/}. Proc.\;London
Math.\;Soc. (1892), {\bf XXIII}, 49--88.
{\em Further note on automorphic functions\/}: 281--295.

\bibitem{308} \mbox{\sc Burnside, W.}
{\em Note on the equation \mbox{$y^2=x\,(x^4-1)$}\/}.
Proc.\;London Math.\;Soc. (1893), {\bf XXIV}, 17--20.

\bibitem{buser} \mbox{\sc Buser, P. \& Silhol, R.}
{\em Geodesics, periods, and equations of real hyperelliptic curves\/}.
Duke Math.\;Journal (2001), {\bf 108}(2), 211--250.

\bibitem{317} \mbox{\sc Costa, A.\,F. \&  Riera, G.}
{\em One parameter families of riemann surfaces of genus two\/}.
Glasgow Math.\;Journal (2001), {\bf 43}, 255--268.

\bibitem{315} \mbox{\sc Chudnovsky, D.\,V. \& Chudnovsky, G.\,V.}
{\em Computer algebra in the service of
mathematical physics and number theory\/}.
Computers in mathematics.
Lect.\;Note in Pure and Appl.\;Math. {\bf 125}, 109--232. Dekker, New-York (1990).

\bibitem{dalaker}\mbox{\sc Dalaker, H.\,H.}
{\em On the Automorphic Functions of the Group $(0,3;2,4,6)$\/}.
Ann.\;Math. (1924), {\bf 25}(3), 241--260.

\bibitem{318} \mbox{\sc Dalzell, D.\,P.}
{\em A note on automorphic functions\/}.
J.\;London Math.\;Soc. (1930), {\bf 5}, 280--282.

\bibitem{dedekind} \mbox{\sc Dedekind, R.} {\em Schreiben an Herrn Borchardt \"uber die
theorie der elliptischen Modulfunctionen\/}.
Borchardt's Journ.\;f\"ur Mathematik (1877), {\bf LXXXIII}, 265--292.

\bibitem{dyck} \mbox{\sc Dyck, W.} {\em Notiz \"uber eine regul\"are
Riemann'sche Fl\"ache vom Geschlechte drei und die zugeh\"orige
,\!,Normalcurve''
vierter\/.} Math.\;Annalen (1880), {\bf XVII}, 510--516.

\bibitem{far} \mbox{\sc Farkas, H.\,M., \ Kopeliovich, Y. \& \ Kra, I.}
{\em Uniformization of modular curves\/.}
Communications in Analysis and Geometry (1996), {\bf 4}(1--2), 207--259.

\bibitem{527} \mbox{\sc Farkas, H.\,M. \&  Kra, I.} {\em Riemann surfaces\/}.
Springer--Verlag (1982).

\bibitem{farkas2000} \mbox{\sc Farkas, H.\,M. \& Kra, I.}
{\em Theta Constants, Riemann Surfaces and the Modular Group\/}.
Graduate Studies  in Math. {\bf 37}, AMS (2001), 531 pp.

\bibitem{551} \mbox{\sc Ford, L.} {\em Automorphic Functions\/}.
McGraw--Hill Book Comp., Inc. New-York (1929).

\bibitem{468} \mbox{\sc Fricke, R.}
{\em Automorphe funktionen mit eischluss der elliptischen
modulfunktionen\/}.
Encyklop\"adie der Mathematischen
Wissenschaften \ldots\ \ (1913), 349--470.
Band II in drei teilen. Analysis.
Verlag und Druck von B.\;G.\;Teubner, Leipzig (1901--1921).

\bibitem{466} \mbox{\sc Fricke, R. \& Klein, F.}
{\em Vorlesungen \"uber die theorie der automorphen Funktionen\/.
{\em {\bf I:}} Die gruppentheoretischen Grundlagen\/}.
Teubner, Leipzig (1897).

\bibitem{467} \mbox{\sc Fricke, R. \&  Klein, F.}
{\em Vorlesungen \"uber die
Theorie der automorphen
Funktionen\/. {\em {\bf II:}} Die Funktionen theorischen Aus f\"uhrungen und
die Anwendungen\/.
Erste Lieferung: Engere Theorie der automorphen Funktionen\/}.
Teubner, Leipzig (1912).

\bibitem{ger} \mbox{\sc Gerstenmeier, C.}
{\em Beitr\"age zur Theorie der
linearen Differentialgleichungen mit $4$ und $5$
singul\"aren Stellen\/.}
Inaugural Dissertation. Erlangen (1910), 88 pp.

\bibitem{479} \mbox{\sc Giraud, G.}
{\em Lecons sur les fonctions automorphes\/}. Paris (1920).

\bibitem{girondo} \mbox{\sc Girondo, E. \& Gonz\'ales-Diez, G.}
{\em On a conjecture of Whittaker concerning uniformization of
hyperelliptic curves\/}. Trans.\;Amer.\;Math.\;Soc. (2004), {\bf 356}(2),
691--702.

\bibitem{golubev} \mbox{\sc Golubev, V.\,V.}
{\em Lectures on analytic theory of differential equations\/}.
Moskva--Leningrad (1950). (in Russian). In German: {\sc Golubew,
W.W.} {\em Vorlesungen der Differetialgleichungen im komplexen\/}.
VEB Deutscher Verlag der Wissenschaften. Berlin (1958).

\bibitem{harnad} \mbox{\sc Harnad, J. \& McKay, J.}
{\em Modular Invariants and Generalized Halphen Systems\/}. SIDE
III---symmetries and integrability of difference equations
(Sabaudia, 1998), 181--195. CRM Proc.\,\,Lecture Notes {\bf 25}, AMS
(2000).

\bibitem{mckay2} \mbox{\sc Harnad, J. \& McKay, J.}
{\em Modular solutions to equations of generalized Halphen type\/}.
Proc.\;Roy.\;Soc.\;London (2000), {\bf 456}(1994), 261--294.

\bibitem{473} \mbox{\sc Harvey, W.\,J.}
{\em Discrete Groups and Automorphic Functions\/}.
Ed. by W.\;J.\;Harvey. Academic Press (1977).

\bibitem{haskell} \mbox{\sc Haskell, M.\,W.}
{\em Ueber die zu der Curve $\lambda^3\mu+\mu^3\nu+\nu^3\lambda=0$
in projectiven Sinne geh\"orende mehrfache Ueberdeckung der
Ebene\/.} Amer.\,Journal of Math. (1890), {\bf 13}(1), 1--52.

\bibitem{342} \mbox{\sc Hejhal, D.\,A.}
{\em Sur les parameteres accesseires pour l'uniformization fuchsienne\/}.
Compt.\;Rend.\;Acad.\;Sci.\;Paris Ser. A--B (1974), {\bf 279},
A695--A697; A713--A716.

\bibitem{343} \mbox{\sc Hejhal, D.\,A.}
{\em Monodromy groups and linear polymorphic functions\/}.
{Acta Mathematica (1975), {\bf 135}(1--2), 1--55.}

\bibitem{hempel} \mbox{\sc Hempel, J.\,A.}
{\em On the uniformization of the $n$-punctured sphere\/}.
Bull.\,London Math.\,Soc. (1988), {\bf 20}(2), 97--115.

\bibitem{hodg}\mbox{\sc Hodgkinson, J.}
{\em Note on the uniformization of hyperelliptic curves\/}.
Journ.\;London Math.\;Soc. (1936), {\bf 11}, 185--192.

\bibitem{411/2}\mbox{\sc  Hutchinson, J.\,I.}
{\em On Certain Relations Among the Theta Constants\/}.
Trans.\;\-Amer.\;\-Math.\; Soc. (1900), {\bf 1}(4), 391--394.

\bibitem{411}\mbox{\sc Hutchinson, J.\,I.}
{\em On a class of automorphic functions.}
Trans.\;Amer.\;Math.\;Soc. (1902), {\bf 3}(1), 1--11.

\bibitem{412}\mbox{\sc Hutchinson, J.\,I.}
{\em On the automorphic functions of the group $(0,3;2,6,6)$\/}.
Trans.\;Amer. Math.\;Soc. (1904), {\bf 5}(4), 447--460.

\bibitem{413}\mbox{\sc Hutchinson, J.\,I.}
{\em A method for constructing the fundamental region of a discontinuous
group of linear transformation\/}.
Trans.\;Amer.\;Math.\;Soc. (1907), {\bf 8}(2), 261--269.


\bibitem{jacobi}  \mbox{\sc Jacobi, C.\;G.\;J.}
{\em Fundamenta Nova Theoriae Functionum Ellipticarum\/.}
Sumtibus Fratrum Borntraeger. Univ. Regiomonti.
K\"onigsberg (1829).
In: {\em Gesammelte Werke\/} {\bf I}, 49--239.
Berlin, Reiner (1882--1891).

\bibitem{jacobiCrelle} \mbox{\sc Jacobi, C.\;G.\;J.}
{\em \"Uber die differentialgleichung, welcher die reihen $1 \pm
2q+2q^4\pm 2q^9+\mathit{etc.},\;\;
2\sqrt[4]{q\mathstrut\,}+2\sqrt[4]{q^9}+2\sqrt[4]{q^{25}}+\mathit{etc}$\/.}
Crelle Journ.\;Reine\;Angew.\;Math. (1847) {\bf 36}, 97--112.

\bibitem{354} \mbox{\sc J\"orgensen, T. \& \ N\"a\"at\"anen, M.}
{\em Surfaces of genus $2$: generic fundamental polygons\/}.
Quart.\;J.\;Math.\;Oxford Ser.\;(2) (1982), {\bf 33}(132), 451--461.

\bibitem{katok} \mbox{\sc Katok, S.}
{\em Fuchsian Groups\/.} The University of Chicago Press (1992).

\bibitem{355} \mbox{\sc Keen, L.}
{\em Canonical polygons for finitely generated Fuchsian groups\/}.
Acta Mathematica (1966), {\bf 115}, 1--16.

\bibitem{klein1} \mbox{\sc Klein, F.}
{\em Ueber die Transformation der
siebenter Ordnung der elliptischen Functionen\/.}
Math.\;Annalen (1878), {\bf XIV}, 428--471.

\bibitem{klein2}  \mbox{\sc Klein, F.}
{\em \"Uber  gewisse Theilwerthe der $\Theta$-Funktionen\/.}
Math.\;Annalen (1881), {\bf XVII}, 565--574.

\bibitem{klein3} \mbox{\sc Klein, F.}
{\em Gesammelte Mathematische Abhandlungen\/} {\bf III}. (1923).
Berlin, Verlag von Julius Springer.

\bibitem{483} \mbox{\sc Klein, F. \&  Fricke, R.}
{\em Vorlesungen \"uber die Theorie
der elliptishen Modulfunktionen.
{\em {\bf I}:} Grundlegung der Theorie;
{\em {\bf II}:} Fortbildung und Anwendung der Theorie\/}.
Teubner, Leipzig (1890--1892).

\bibitem{knopp2} \mbox{\sc Knopp, M.\,I.}
{\em On Abelian integrals of the second kind and modular forms\/}.
Amer.\;J.\;Math. (1962), {\bf 84}(4), 615--628.

\bibitem{knopp} \mbox{\sc Knopp, M.\,I.}
{\em Modular Functions in Analytic Number Theory\/}. Chelsea Pbl.,
New York (1993).

\bibitem{367} \mbox{\sc Kuusalo, T. \&  N\"a\"at\"anen, M.}
{\em Geometric uniformization in genus $2$\/}.
Annales Acad.\;Sci.\; Fennic\ae\ Ser.\;A.\;I.\;Math.\;(1995),
{\bf 20}, 401--418.

\bibitem{kuyk} \mbox{\sc Kuyk, W.} (Ed.) {\em Modular Functions of One
Variable I\/}. Lect.\,Notes in Math. {\bf 320} (1973).

\bibitem{lehner} \mbox{\sc Lehner, J.}
{\em Discontinuous groups and automorphic functions\/}.
Providence: Amer.\;Math.\;Soc. (1964).

\bibitem{LS} \mbox{\sc Lyndon, R. \& Schupp, P.}
{\em Combinatorial group theory\/}.
Springer--Verlag (1977).

\bibitem{magnus} \mbox{\sc Magnus, W., Karrass, A. \& Solitar, D.}
{\em Combinatorial group theory\/}. New York: Wiley (1966).

\bibitem{mckay} \mbox{\sc McKay, J. \& Sebbar, A.}
{\em Fuchsian groups, automorphic functions and Schwarzians\/}.
Math.\;Annalen (2000), {\bf 318}(2), 255--275.

\bibitem{mckean} \mbox{\sc McKean, H. \& Moll, V.}
{\em Elliptic curves\/.} Cambridge University Press (1997).

\bibitem{447} \mbox{\sc Morris, R.}
{\em On the automorphic functions of the group $(0,3;l_1, l_2,
l_3)$\/}. Trans.\;Amer.\;Math.\; Soc. (1906), {\bf 7}(3),
425--448.

\bibitem{mumford} \mbox{\sc Mumford, D.}
{\em The Red Book of Varieties and Schemes\/}. Lect.\,Notes in Math.
{\bf 1358}. 2-nd ed. Springer (1999).

\bibitem{378} \mbox{\sc Mursi, M.}
{\em On the uniformization of algebraic curves of genus $3$\/}.
Proc.\;Edinburgh Math.\;Soc. (2) (1930), {\bf 2}, 101--107.

\bibitem{379} \mbox{\sc Mursi, M.}
{\em A note on automorphic functions\/}.
J.\;London Math.\;Soc. (1931), {\bf 6}, 166--169.

\bibitem{539} \mbox{\sc Nevanlinna, R.} {\em Uniformisierung\/}.
Springer--Verlag (1953).

\bibitem{PL} \mbox{\sc Plemelj, J.}
{\em Problems in the sense of Riemann and Klein.\/}
Interscience Publishers (1964).

\bibitem{poincare1} \mbox{\sc Poincar\'e, H.}
{\em Sur les fonctions Fuchsiennes\/.} Compt.\;Rend.\;Acad.\;Sci.
(1881), {\bf 92}, 333--335, 395--398, 1198--1200;
{\em Sur une
nouvelle application et quelques propri\'et\'es importantes des
fonctions Fuchsiennes\/.} Compt.\;Rend.\;Acad.\;Sci. (1881), {\bf
92}, 859--861.



\bibitem{poincare3} \mbox{\sc Poincar\'e, H.}
{\em Sur les fonctions Fuchsiennes\/.}
Compt.\;Rend.\;Acad.\;Sci. (1881), {\bf 92}, 1484--1487.

\bibitem{poincare4} \mbox{\sc Poincar\'e, H.}
{\em Sur une fonctions analogue aux fonctions modulaires\/.}
Compt.\;Rend. Acad. Sci. (1881), {\bf 93},  138--140.

\bibitem{P82} \mbox{\sc Poincar\'e, H.}
{\em Th\'eorie des groupes fuchsiennes\/}.
Acta Mathematica (1882), {\bf 1}, 1--61.

\bibitem{P84} \mbox{\sc Poincar\'e, H.}
{\em Sur les groupes des \'equations lin\'eaires\/}.
Acta Mathematica (1884), {\bf 4}, 201--311.

\bibitem{pl} \mbox{\sc Poincar\'e, H.} {\em Les fonctions Fuchsiennes
et l'\'equation $\Delta u=e^u$\/}. Journal de Math\'ematiques
(1898), s\'er. 5, {\bf 4}, 137--230.

\bibitem{390} \mbox{\sc Rankin, R.\,A.}
{\em On horocyclic groups\/}.
Proc.\;London Math.\;Soc. (3) (1954), {\bf IV}, 219--234.

\bibitem{393} \mbox{\sc Rankin, R.\,A.}
{\em The Differential Equations Associated with the Uniformization
of Certain Algebraic Curves\/}. Proc.\;Royal Soc.\;Edinburgh.
Section A (1958), {\bf LXV}, 35--62.

\bibitem{395} \mbox{\sc Rankin, R.\,A.}
{\em The Schwarzian derivative and uniformization\/}.
J.\;Analyse Math. (1958), {\bf 6}, 149--167.

\bibitem{494} \mbox{\sc Rankin, R.\,A.} {\em Modular forms and functions\/}.
Cambridge Univ.\,Press (1977). London--New-York--Melbourne.

\bibitem{397} \mbox{\sc Rankin, R.\,A.}
{\em Burnside's uniformization\/}.
Acta Arithmetica (1997), {\bf 79}(1), 53--57.

\bibitem{riemann} \mbox{\sc Riemann, B.}
{\em Bernhardt Riemann's Gesammelte Mathematische Werke und
Wissenschaftliche Nachlass\/,} ed. R.\;Dedekind and H.\;Weber
(1892) with {\em Nachtr\"age\/,} ed. M.\;Noeter and W.\;Wirtinger.
Leipzig, B.\;J. Teubner, (1902).

\bibitem{496} \mbox{\sc Schlesinger, L.}
{\em Automorphe Funktionen\/}.
Walter de Gruyter \& Co. Berlin w.\,10 und Leipzig (1924).

\bibitem{405} \mbox{\sc Schottky, F.}
{\em Ueber die conforme Abbildung mehrfach zusammenh\"angender
ebner Fl\"achen\/.} Inaugural--Dissertation. Berlin (1875), 61 pp.
Crelle Journ.\;Reine Angew.\;Math. (1877),
{\bf LXXXIII}, 300--351.

\bibitem{smirnov} \mbox{\sc Smirnov, V.\,I.}
{\em The inversion problem of linear differential equation
of second order with four singular points\/}. Dissertation.
Petrograd (1918), 307 pp. (in Russian)

\bibitem{takh1} \mbox{\sc Takhtajan, L.\,A.}
{\em Uniformization, Local Index Theorem, and Geometry of the Moduli Spaces
of Riemann Surfaces and Vector Bundles\/}.
Proc.\;Sympos.\;Pure Math. {\bf 49}, part 1.
Amer.\;Math.\;Soc., Providence, R1 (1989), 581--596.

\bibitem{takh2} \mbox{\sc Takhtajan, L.\,A.}
{\em A simple example of Modular forms as $\tau$-func\-ti\-ons
for integrable equations\/}.
Theor.\;Math.\;Phys. (1992), {\bf 93}(2), 330--341.

\bibitem{tah3} \mbox{\sc Takhtajan, L. \& Zograf, P.}
{\em Hyperbolic $2$-spheres with conical singularities, accessory
parameters and K\"ahler metrics on $\mathcal{M}_{0,n}$\/.}
Trans.\,Amer.\,Math.\,Soc. (2002), {\bf 355}(5), 1857--1867.

\bibitem{tannery} \mbox{\sc Tannery, J. \& Molk, J.}
{\em Elements de la theorie des fonctions elliptiques\/.}
{\bf I--IV}. Paris, Gauithier--Villars (1893--1902).

\bibitem{thome} \mbox{\sc Thom\'e, L.\,W.}
{\em Ueber lineare Differentialgleichungen
mit mehrwerthigen algebraischen Coefficienten\/.} 
Journ.\;Reine\;Angew.\;Math. (1895), {\bf CXV}, 33--52, 119--149;
(1898), {\bf CXIX}, 131--147.

\bibitem{417} \mbox{\sc Weber, H.}
{\em Ein Beitrag zu Poincar\'e's Theorie der Fuchs'schen Functionen\/}.
G\"ottingen Nachrichten (1886), 359--370.

\bibitem{weber} \mbox{\sc Weber, H.}
{\em Lehrbuch der Algebra\/.} {\bf I--III}, 3-rd ed. Chelsea, New-York.

\bibitem{weyl} \mbox{\sc Weyl, H.} {\em The concept of a Riemann surface\/}.
3-rd ed. Addison--Wesley Publ., Inc. (1955).

\bibitem{418} \mbox{\sc Whittaker, E.\,T.}
{\em On the Connexion of Algebraic Functions with Automorphic Functions\/}.
Phil.\;Trans. of the Royal Soc. of London (1898), {\bf 192A}, 1--32.

\bibitem{420} \mbox{\sc Whittaker, E.\,T.}
{\em On hyperlemniscate functions, a family of automorphic functions.}
J.\;London Math.\;Soc. (1929), {\bf 4}, 274--278.

\bibitem{watson} \mbox{\sc Whittaker, E.\,T. \& Watson, G.\,N.}
{\em A course of modern analysis\/.} 4-th ed. Cambridge Univ.\;Press (1927).

\bibitem{421} \mbox{\sc Whittaker, J.\,M.}
{\em The uniformisation of algebraic curves\/}.
J.\;London Math.\;Soc. (1930), {\bf 5}, 150--154.

\bibitem{young1} \mbox{\sc Young, J.\,W.}
{\em On the group of sign $(0,3;2,4,\infty)$ and the functions
belonging to it\/}. Trans. Amer.\;Math.\;Soc. (1904), {\bf 5}(1),
81--104.

\bibitem{young2}\mbox{\sc Young, J.\,W.}
{\em On class of discontinuous  $\zeta$--groups defined by normal curves of
the fourth order in a space of four dimensions\/}.
Rend.\;Circ.\;Matem.\;Palermo (1907), {\bf XXIII}, 97--106.

\bibitem{tah5} \mbox{\sc Zograf, P.\,G. \& Takhtadzhan, L.\,A.}
{\em Action of the Liouville equation generation function
for accessory parameters and the potential of the Weil--Petersson metric
on Teichm\"uller space\/}.
Func.\;Anal.\;Appl. (1985), {\bf 19}(3), 67--68.

\bibitem{tah4} \mbox{\sc Zograf, P.\,G. \& Takhtajan, L.\,A.}
{\em On Liouville's equation, accessory parameters, and the geometry of
Teichm\"uller spaces for riemann surfaces of genus $0$\/}.
Math.\;USSR Sb. (1987), {\bf 60}(1), 143--161.



\end{thebibliography}
\end{document}